\documentclass[preprint]{article}

\usepackage[utf8]{inputenc}
\usepackage[english]{babel}
\usepackage[T1]{fontenc} % Font encoding
\usepackage{graphicx}
\usepackage{eso-pic} % For the background picture on the title page
\usepackage{subfig} % Numbered and caption subfigures using \subfloat
\usepackage[font=small]{caption} % Coloured captions
\usepackage{transparent}
\usepackage{amsmath}
\usepackage{amssymb}
\usepackage{amsthm}
\usepackage{amsopn}
\usepackage{bm}
\usepackage[overload]{empheq}  % For braced-style systems of equations
\usepackage{booktabs}
\usepackage{siunitx} 
\usepackage{cancel}
\usepackage{makecell}
\usepackage{pdfpages}

\newcommand \diiv {\nabla \cdot}
\DeclareMathOperator{\tr}{tr}

\newcommand \Ra {\mathit{Ra}}
\newcommand \Sh {\mathit{Sh}}
\newcommand \Kb {k_b/(\phi\mu)\;}
\newcommand \Kf {k_f/(\phi\mu)\;}
\newcommand \Kd {k_d/(\phi\mu)\;}
\newcommand{\dd}{\,d}
\newcommand{\de}{\,\partial}
\newcommand{\noe}{\mathit{NoE}}
\newcommand{\hd}[1]{\makecell[c]{ #1 }}
\newtheorem{remark}{Remark}
\providecommand{\keywords}[1]{\textbf{\textit{Keywords:}} #1}

\begin{document}

\author{Arash Andrea Roknian$^1$ \and Anna Scotti$^2$ \and Alessio Fumagalli $3$}
\date{MOX - Dipartimento di Matematica ``F. Brioschi'', Politecnico di
Milano, via Bonardi 9, 20133 Milan, Italy\\
    $^1$ arashandrea.roknian@mail.polimi.it\\%
    $^2$ anna.scotti@polimi.it\\
    $^3$ alessio.fumagalli@polimi.it
    \today
}

\title{Free convection in fractured porous media: a numerical study}

\maketitle

\begin{abstract}
The objective of this study is to better understand the influence of fractures
on the possibility of free convection in porous media. To this aim, we introduce a mathematical model for density driven flow in the presence of fractures, and the corresponding numerical approximation. In addition to the direct numerical solution of the problem we propose and implement a novel method for the assessment of convective stability through the eigenvalue analysis of the linearized numerical problem. The new method is shown to be in agreement with existing literature cases both in simple and complex fracture configurations. With respect to direct simulation in time, the results of the eigenvalue method lack information about the strength of convection and the steady state solution, they however provide detailed (quantitative) information about the behaviour of the solution near the initial equilibrium condition. Furthermore, not having to solve a time-dependent problem makes the method computationally very efficient. Finally, the question of how the porous matrix interacts with the fracture network to enable free convection is examined: the porous matrix is shown to be of key importance in enabling
convection for complex fracture networks, making stability criteria based on the fracture network alone somewhat limited in applicability.
\end{abstract}



\section{Introduction}
\label{intro}

The study of flow in porous media finds application in many different areas ranging from industrial, to biomedical and environmental applications. 
In the context of geological porous media many relevant problems such as geothermal engineering and contaminant transport are nowadays addressed with the help of numerical simulations, due to the necessity of estimating flow rates and pathways.

In many cases, due to either temperature or the presence of solutes, fluids can present density variations in space and time.
Temperature gradients are particularly important in geothermal applications, 
while solute concentration gradients are found for instance in the study of contaminant plume migration or seawater intrusion phenomena.

In situations of inverse density gradients (denser fluid on top) the fluid may spontaneously develop unstable convective plumes.
It is particularly important to predict the possible onset of such instabilities due to the strong rates at which convection can transport solutes with respect to diffusion alone.

The first studies dealing with free convection in porous media \cite{lapwood_1948,hortonrogers} examined simple scenarios of inverse density gradients. The possibility of convection is linked to the Rayleigh number, which, accounting for both fluid and porous medium parameters, quantifies the strength of convection against that of diffusion.
For large values of the Rayleigh number, there is the possibility of free convection.

In recent years, different authors have tried to understand how the presence of fractures influences the possibility and strength of free convection \cite{shikaze98,shafabakhsh19,vg14,vg15,nield}.
While high-density fracture networks can be mostly treated through averaging/upscaling of porous media properties, low-density fracture networks can manifest unexpectedly high convection strength due to the particular geometry of the network.
Even for simple geometries, e.g. for regular horizontal or vertical fracture grids, the unstable nature of convective plumes makes prediction very hard if not impossible \cite{shikaze98}.

The main focus of this study is the development of a new method based on the eigenvalue analysis of equilibrium solutions for assessing the possibility of convection in fractured porous media. Leveraging this new method will enable us to better understand the particular mechanisms by which fractures enhance convection. This idea is at the core of the results provided by \cite{lapwood_1948}, where the analysis is carried out analytically for homogenoues media, and in \cite{rees} for layered porous media. Here, given the geometrical complexity of the domain, eigenvalues will be computed by a suitable numerical algorithm starting from the discretized problem.

We start in section \ref{sec:model} with a detailed description of the mathematical model both in homogeneous and fractured media.
As the focus of this study is the effect of fracture geometry, simple constitutive models are preferred over more elaborate (and more accurate) models.
The Darcy law is used as constitutive law for filtration, Fick's law is used to model solute diffusion and density is modeled as a linear function of solute concentration.

Particular care is given to the averaging procedure used for the dimensional reduction of fractures, consistent with the mixed-dimensional approach presented in \cite{martin05}.

In section \ref{sec:discretization} we present the discretization of the continuous model using the finite volume method in space and the implicit Euler scheme in time.
The particular finite volume method used in the implementation is the Multipoint Flux Approximation method \cite{aavatsmark02} for both the fluid mass conservation and the solute transport problems. 
The MPFA method is particularly appropriate due to its mass conservation properties and consistency for general grids. 
Again, density dependence must be treated with care for the numerical implementation to be consistent \cite{starnoni19}.

Two methods will be described for assessing the possibility of convective cell formation. 

Section \ref{sec:direct} describes what we call the direct method of assessing stability: starting from a non-equilibrium initial condition and integrating the time-dependent flow equations until steady state may result in either convective motion or a diffusive equilibrium solution. The two regimes are easily distinguished by measuring the amount of solute transport through the domain.

Section \ref{sec:eig} describes the eigenvalue method for assessing stability. By linearizing the problem and numerically studying the eigenvalues of the resulting discrete system, we can assess the stability of arbitrary perturbations to any equilibrium solution without the need of solving the time-dependent problem.

Sections \ref{sec:exp-elder} and \ref{sec:exp-hrl} test the two methods respectively with the Elder and HRL problems relying on the results of \cite{vg14} and \cite{diersch02}: both are well known benchmark cases for density driven flow.
Section \ref{sec:exp-hrl-3d} validates and examines a three-dimensional generalization of the HRL problem, based on the studies in \cite{nield, vg15}.
In the HRL scenarios, the complementarity of the two approaches described above will enable us to better understand and examine the peculiarities of free convection in the presence of fractures. 

The variety of scenarios is also useful for understanding the advantages and limitations of the eigenvalue method with respect to the direct method. A comparison of their computational cost is presented in section \ref{sec:computational}.

Finally, section \ref{sec:conclusioni} will conclude this study with a critical evaluation of what it has accomplished and suggest different options for further developments.

\section{Mathematical model}\label{sec:model}

The problem of our interest stems from the coupling between solute transport and density driven flow in a porous medium.
The mathematical model is based on the one described in the review paper by Diersch and Kolditz \cite{diersch02}, which will be here extended to account for the presence of fractures in the domain, since our main goal is to understand the impact of fracture networks on the onset of free convection.

Let us begin by considering a generic advection-diffusion equation expressing solute mass conservation in a porous medium:
\begin{equation}
\label{eq:homogenous_solute_conservation}
\partial_t (\rho \phi \omega) + \nabla \cdot(\rho \phi \omega \bm{u} )+\nabla \cdot(\rho \phi \bm{i} )=0 ,
\end{equation}

\noindent where $\rho\;[\mathrm{M}/\mathrm{L}_v^3]$ is the density of the fluid in which the solute is dissolved, $\omega\;[\mathrm{M_s}/\mathrm{M}]$ is the concentration of the solute, $\bm {u} \; [\mathrm{L}/\mathrm{T}]$ is the flow velocity, $\rho \phi \bm{i} \; [\frac{\mathrm{M_s}}{\mathrm{L}^2 \mathrm{T}}]$ is the solute mass flux due to diffusion and $\phi \; [\mathrm{L_v^3}/\mathrm{L^3}] $ is the porosity of the medium. 
The subscript $v$ in the length dimension $L$ is used to differentiate void volumes from total volumes, while $s$ is used to differentiate solute mass from fluid mass.

For the diffusive flux a Fick-type law is used: $ \bm{i}=  -D \nabla \omega . $
The diffusion tensor $D \;[\mathrm{L}^2/\mathrm{T}]$ is usually modeled as the sum of two components: 
a part related to molecular diffusion and a part related to mechanical dispersion (function of the fluid velocity): $ D= D_d (\bm{u})+ D_m I$.
In this study only the part related to molecular diffusion has been considered, thus 
\begin{equation}
\label{eq:homogeneous_diffusion}
\bm{i}=-D_m \nabla \omega .
\end{equation}
In the following, it will be useful to indicate the combined advective and diffusive solute flux as
$\bm{q} \;[\frac{\mathrm{M_s}}{\mathrm{L}^2 \mathrm{T}}]$: 
\begin{equation}
\label{eq:dfn_q}
\bm{q} = \omega \bm{u} + \bm{i} .
\end{equation}
The boundary conditions for equation (\ref{eq:homogenous_solute_conservation}) can be 
of Dirichlet type where we set the value of the solute concentration $ \omega = \omega_D $, 
or of Neumann type where the total mass flux in the normal direction is specified $ \bm{q} \cdot \bm{n} = q_N $, where $\bm{n}$ is the unit normal (by convention pointing outwards from the domain).

To close the system we still need an expression for the fluid velocity and a constitutive equation for density as a function of primary variables.
The first can be determined as the solution of the Darcy problem, which is the model used to describe filtration in porous media. In particular we have a mass conservation law for the fluid:
\begin{equation}
\label{eq:homogenous_mass_conservation}
\partial_t (\rho \phi)+\nabla \cdot(\rho \phi \bm{u} )=0,
\end{equation}
and a constitutive law classically used in the context of porous media, the Darcy law:
\begin{equation}
\label{eq:homogenous_darcy}
\bm{u} = \frac{k}{\phi \mu} (- \nabla p+\rho \bm{g}) .
\end{equation}

where $\bm{g}=-g\bm{e}_z$ is the gravity acceleration vector, $k$ is the permeability tensor of the porous medium $k \;[\mathrm{L}_v^2]$, and $\mu \;[\mathrm{M}/\mathrm{L}_v\mathrm{T}]$ is the viscosity of the fluid.

The boundary condition for (\ref{eq:homogenous_mass_conservation}) can be 
of Dirichlet type where the pressure is specified: $ p = p_D $ 
or of Neumann type where the normal flow velocity is prescribed, $ \bm{u} \cdot \bm{n} = u_N $.

For the density constitutive law, we will assume a linear dependence on the concentration $\omega$
\begin{equation}
\label{eq:density}
\rho = \rho(p, \omega) = \rho_0 (1 + \alpha \omega) ,
\end{equation}
where $\rho_0$ is a reference value and $\alpha$ a dilation coefficient.

A large part of the following discussion applies unchanged to flows driven by temperature gradients instead of solute concentration gradients.
Indeed, apart from a reinterpretation of some of the quantities introduced above e.g. thermal conductivity instead of diffusivity, 
the same model can be applied.

\subsection{Oberbeck-Boussinesq approximation}

The analysis and the solution of system (\ref{eq:homogenous_solute_conservation}) -(\ref{eq:density})  can be substantially simplified by assuming the solution satisfies what is known as the Oberbeck-Boussinesq (OB) approximation.

The OB approximation consists in neglecting all but the most important among the different density related nonlinearities 
in the system.
In particular we will take $\rho = \rho_0 (1 + \alpha \omega) \approx \rho_0 $ everywhere except in the gravity term. 

Using this approximation, and after some algebraic simplifications the system reduces to 
\begin{equation}
\label{eq:homogeneous}
\left\{
\begin{aligned}
&\nabla \cdot \bm{u} = 0  && \text{in } \Omega, \\
&\partial_t \omega + \nabla \cdot \bm{q} = 0 & \quad & \text{in } \Omega, \\
&\bm{u}= \frac {k}{\phi \mu} (- \nabla p_e + \rho_0 \alpha \omega \bm{g}) && \text{in } \Omega, \\
&\bm{q}= - D \nabla \omega + \omega \bm{u} && \text{in } \Omega, \\
\end{aligned}
\right.
\end{equation}
with boundary conditions
\begin{equation*}
\left\{
\begin{aligned}
    & \bm{u} \cdot \bm{n} = u_N & \quad & \text{on } \partial\Omega_N^p, \\
    & p = p_D && \text{on } \partial\Omega_D^p,
\end{aligned}
\right.\qquad\left\{
\begin{aligned}
    & \bm{q} \cdot \bm{n} = q_N && \text{on } \partial\Omega_N^\omega, \\
    & \omega = \omega_D && \text{on } \partial\Omega_D^\omega, \\
\end{aligned}
\right.
\end{equation*}
where $\partial\Omega = \partial\Omega_N^p \cup \partial\Omega_D^p $, $\partial\Omega = \partial\Omega_N^\omega \cup \partial\Omega_D^\omega$.

Among the simplifications to reach the reduced system, we have replaced the unknown $p$ with the excess pressure 
$ p_e $ by removing the hydrostatic component $p_h= \rho_0 g (y_0 - y)$ (defined for the reference density, in the absence of solute):
\begin{equation*}
\begin{aligned}
    p  = p_e + p_h           &= p_e + \rho_0 g (y_0 - y) , \\
    -\nabla p + \rho \bm{g}  &= -\nabla p_e - \cancel{\rho_0 \bm{g}} +\cancel{\rho_0 \bm{g}} + \rho_0 \alpha \omega \bm{g}   
\end{aligned}
\end{equation*}    

For readability, in what follows we will drop the subscript $e$ and use the variable $p$ to indicate the excess pressure $p_e$.

\subsection{Horton-Rogers-Lapwood (HRL) problem}
\label{sec:hrl}

\begin{figure}
    \centering
    \includegraphics[width=\textwidth]{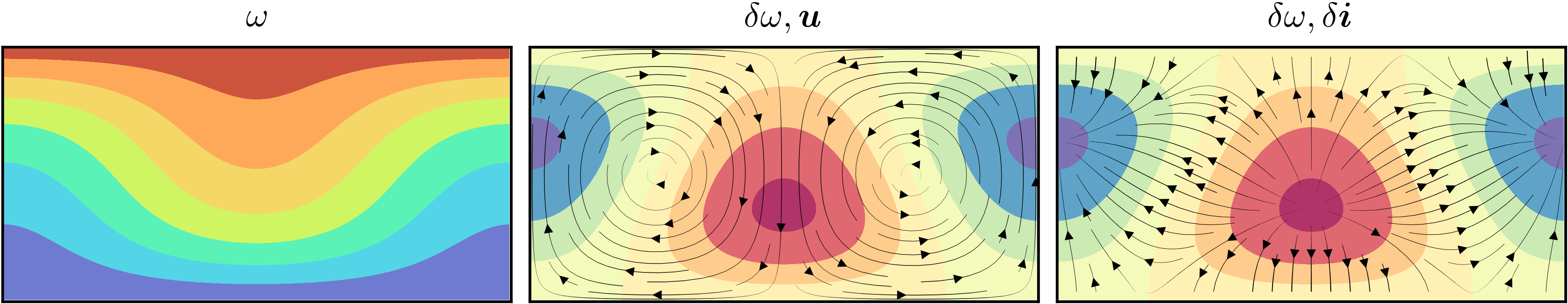}
    \caption{Convective motion. On the left, concentration profile. At center, perturbation of concentration and flow velocity. On the right, the diffusive flow trying to restore the concentration imbalance. }
    \label{fig:ra}
\end{figure}

The HRL problem is a simple scenario aiming to study the possible onset and strength of natural convection.
The idea is to impose, through the boundary conditions, a layer of heavier fluid overlaid on lighter fluid in a two-dimensional vertical cross section of a homogeneous porous medium.
The density gradient can be caused by temperature difference (as in the original description \cite{lapwood_1948}), or by solute concentration difference.
In a situation of inverse mass gradient, the fluid will form convective cells only if diffusivity is small enough to allow it.
The contrast between convection and diffusion speed, respectively $v_g$ and $v_d$, is described by the Rayleigh number
\begin{equation}
\label{eq:rayleigh}
\Ra = \frac{v_g}{v_d} = \frac{ \frac{k}{\phi \mu} \rho_0 \alpha \omega_{max} g } {\frac {D}{H}} ,
\end{equation}

where $H$ is the height of the domain, and other quantities have been introduced in the previous section.
In the original presentation of this problem \cite{lapwood_1948}, the possibility of convective motion is linked to the value of $\Ra$: $\Ra_c = 4 \pi^2$ is the critical number i.e. a threshold such that, for $ Ra > Ra_c $ convection occurs. 
Why is free convection only possible for such values of the Rayleigh number?
For high enough diffusion, convective cells cannot sustain the concentration/temperature difference and thus they simply decay.
Figure \ref{fig:ra} illustrates how the diffusive flow tries to restore concentration imbalances caused by the convective motion.

Note that convective motion in general enhances solute transport. 
For this reason, an indicator of the presence of convection is the Sherwood number $\Sh$, defined as
\begin{equation}
\label{eq:sherwood}
\Sh= \frac{\int_A \bm{i} \cdot \bm {n}}{\int_A \bm{i}_0 \cdot \bm {n}} =\frac{\int_A \bm{i} \cdot \bm {n}}{\frac {D}{H} A} ,
\end{equation}
where $A$ is the diffusive inflow surface (top boundary) and $\bm{i}_0$ is the diffusive flow in the absence of convection. When convection is not present inside the domain, $ \bm{i} = \bm{i_0} $ and thus $\Sh=1$. Note that the boundary conditions for the HRL problem being of zero fluid mass flow all around, the solute transport on the inflow surface $A$ is entirely due to diffusion, regardless of whether convection is present inside the domain. However, when convection is present, the (overall) larger concentration gradient at $A$ makes $ \bm{i} > \bm{i_0} $ and $ \Sh > 1 $.

\subsection{Dimensional reduction}

Porous media often present heterogeneities in their material properties.
A particularly strong kind of heterogeneity are fractures: regions of different material properties, with negligible aperture (thickness) with respect to both their length and the characteristic lengths of the medium.
Very different material properties in fractures, such as permeability, can compensate for their small dimensions so that, overall, fractures can strongly influence the behaviour of flow in the medium. Let us define 

$$
k = k(\bm{x}) = \left\{
\begin{aligned}
    & k_b \enspace&& \bm{x} \in \Omega_b, \\
    & k_f         && \bm{x} \in \Omega_f,
\end{aligned}
  % \begin{array}{ll}
  % \end{array}
\right. 
$$
where $\Omega_f$ denotes the fracture, $\Omega_b$ the bulk medium and $ \Omega = \Omega_b \cup \Omega_f $.

% SIDE CAPTION FIGURA GRANDE
\begin{figure}
\centering
\includegraphics[width=0.5\textwidth]{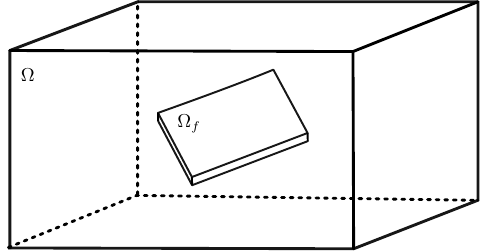}
\caption{Domain $\Omega$ containing a single planar fracture $\Omega_f$.}
\label{fig:dominio-fratturato} 
\end{figure}

The approach known as discrete fracture modeling aims to solve this problem by treating fractures as separate domains, often of reduced dimensionality, coupled to the bulk medium by coupling conditions at their interface.

\begin{figure}
\centering
\includegraphics[width=0.8\textwidth]{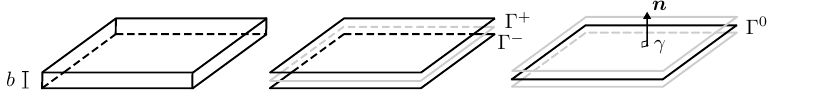}
\caption{$\Omega_f$ is a rectangular fracture of height $b$. $\Gamma^+$ and $\Gamma^-$ are both subsets of its boundary: $\Gamma^\pm \subset \partial \Omega_f$. }
\label{fig:domini-interfacce}
\end{figure}

We will now rewrite the fluid mass conservation equation and associated Darcy law as a mixed-dimensional equation. 
To keep our derivation simple, we will treat a single two-dimensional fracture as illustrated in figure \ref{fig:dominio-fratturato}. The domain $\Omega$ is split between bulk medium and fracture $\Omega = \Omega_b \cup \Omega_f \cup \Gamma $, where $ \Gamma = \overline{\Omega_b} \cap \overline{\Omega_f} $ denotes the interface between subdomains. We assume that $\Omega_f$ can be expressed as
\begin{align*}
& \Omega_f = \{\bm{x} \in \Omega:\bm{x}=\bm{\gamma} + \alpha \bm{n},\: \gamma \in \Gamma^0,\: \alpha \in ( -b/2, b/2 ) \} , \\
& \Gamma^0 = \{\bm{x} \in \Omega:\bm{x}=\bm{x}_0 + s_1(\bm{x}_1 - \bm{x}_0) + s_2(\bm{x}_2 - \bm{x}_0),\: s_i \in (0,1) \}  ,
\end{align*}
where we have assumed a planar midsurface $\Gamma^0$ as shown in figure \ref{fig:domini-interfacce}. The lateral boundary of $\Omega_f$ will be treated differently from the top and bottom boundaries. We define
\begin{align*}
& \Gamma^\pm = \{ \bm{x} \in \partial\Omega_f : \bm{x} = \bm{\gamma} \pm \frac{b}{2} \bm{n},\: \gamma \in \Gamma^0 \} , \\
& \Sigma     = \{ \bm{x} \in \partial\Omega_f : \bm{x} = \bm{\sigma} + \alpha \bm{n}, \: \sigma \in \partial\Gamma^0, \: \alpha \in ( -b/2, b/2 ) \}  ,
\end{align*}
such that $\Gamma = \Gamma^\pm \cup \Sigma $. 
Superscripts + and - will also be used to denote quantities evaluated on $\Gamma^\pm$.

Starting from the continuity equation and the associated Darcy law in the two domains, we enforce pressure continuity and conservation of mass across the interface $ \Gamma $:

\begin{equation}
\label{eq:hetero_darcy}
\left\{
\begin{aligned}
        &\diiv \bm{u}_b = 0 && \text {in } \Omega_b ,\\
        &\bm{u}_b= \Kb (- \nabla p_b + \rho_0 \alpha \omega_b \bm{g}) & \qquad & \text {in } \Omega_b ,\\[1mm]
        &\diiv \bm{u_f} = 0 && \text {in } \Omega_f ,\\
        &\bm{u_f}= \Kf (- \nabla p_f + \rho_0 \alpha \omega_f \bm{g}) && \text {in } \Omega_f ,\\[1mm]
        &\tr p_b = \tr p_f && \text {on } \Gamma ,\\
        &\tr \bm{u}_b \cdot \bm{n}_b = \tr \bm{u_f} \cdot \bm{n}_b && \text {on } \Gamma ,\\
\end{aligned}
\right.
\end{equation}
where $\bm{n}_b$ is the normal vector defined on $\Gamma$, exiting from $\Omega_b$.
At the end of dimensional reduction, the equations solved on the three-dimensional domain $\Omega_f$ will be replaced by the solution of (different) equations on the two-dimensional domain $\Gamma^0$. 
For the sake of simplicity, we will suppose that fractures differ from the bulk only in their permeability: all the other material parameters will be common to both. 

We proceed by integrating the mass conservation equation in the fracture across the aperture.
$$
0 = \int_{-b/2}^{b/2} \diiv \bm{u_f} = \int_{-b/2}^{b/2} \diiv (\mathrm{T} \bm{u_f} + \mathrm{N} \bm{u_f}) \\
  = \nabla_{\tau} \cdot \bm{u_\gamma} - (\bm{u}_b \cdot \bm{n}_b)^+ - (\bm{u}_b \cdot \bm{n}_b)^- ,
$$
where we have denoted by $\mathrm{T}$, $\mathrm{N}$ the projection operators specific to the fracture: 
$ \mathrm{N} \bm{v} = (\bm{v} \cdot \bm{n}) \bm{n} $ and 
$ \mathrm{T} \bm{v} = (\mathrm{I} - \mathrm{N}) \bm{v} $, where $\bm{n}$ is the normal vector to the fracture plane (as indicated in figure \ref{fig:domini-interfacce}), and $\mathrm{I}$ is the identity operator.
We also denoted the in-plane gradient operator by $\nabla_\tau = \mathrm{T}\nabla$, and by $\bm{u}_\gamma$ the integral of the tangential flow field:
$$\bm{u}_\gamma = \int_{-b/2}^{b/2} \mathrm{T} \bm{u_f}, \quad [\mathrm{L}^2/\mathrm{T}]. $$
Indeed, given these definitions,
$$ 
\int_{-b/2}^{b/2} \diiv (\mathrm{T} \bm{u_f}) = 
\int_{-b/2}^{b/2} \nabla_{\tau}\cdot (\mathrm{T} \bm{u_f}) = 
\nabla_{\tau}\cdot \int_{-b/2}^{b/2} \mathrm{T} \bm{u_f} =
\nabla_{\tau}\cdot \bm{u}_\gamma.
$$

% By integrating across the aperture, the balance of fluxes has split into an in-plane contribution and a contribution from the surrounding porous medium. The in-plane part can be written as the in-plane divergence $\nabla_\tau \cdot $ of a vector field $\bm{u}_\gamma $. The divergence however now only operates inside the fracture plane and, since $ \nabla_\tau = \mathrm{T} \nabla $, $ \bm{u}_\gamma$ is the integral of the tangential flow field: 

To obtain a law for $\bm{u}_\gamma$, we integrate the in-plane component of the Darcy law:
\begin{align*}
\bm{u_\gamma} 
= \int_{-b/2}^{b/2} \mathrm{T} \bm{u_f} 
= \frac{k_f}{\phi \mu} \int_{-b/2}^{b/2} \mathrm{T} (- \nabla p_f + \rho_0 \alpha \omega_f \bm{g})
= \frac{b k_f}{\phi \mu} (-\nabla_{\tau} p_\gamma + \rho_0 \alpha \omega_\gamma \bm{g_\tau}) ,
\end{align*}
where $p_\gamma = \frac{1}{b} \int_{-b/2}^{b/2} p_f $, $\omega_\gamma = \frac{1}{b} \int_{-b/2}^{b/2} \omega_f $ and $ \bm{g_\tau} = \mathrm{T} \bm{g}$. Differently from $ \bm{u}_\gamma $, scalar quantities in the fracture are averaged across the fractures thus maintaining the same dimensions as the corresponding bulk quantities.
Coupling conditions at the interface must also be expressed in terms of the averaged variables:
\begin{equation*}
(\tr \bm{u}_b \cdot \bm{n}_b)^\pm
 \approx \frac {b k_f}{\phi \mu} (\frac{p_b^\pm - p_\gamma}{b/2} + \rho_0 \alpha \omega_b^\pm g_n).
\end{equation*}

We conclude by replacing $\Omega_f$ with its center plane $\Gamma^0$ and extending $\Omega_b$:  $\Omega = \Omega_b \cup \Omega_f \cup \Gamma \approx (\Omega \setminus \Gamma^0) \cup \Gamma^0 \cup \Gamma^\pm$. Even though the fracture domain $ \Omega_f $ and the interface $\Gamma$ have collapsed on one another geometrically, the two play distinct roles and must be kept conceptually separate.

Collecting the last steps, we can write the Darcy part of the mixed-dimensional system:
\begin{equation}
\label{eq:mixed_darcy}
\left\{
\begin{aligned}
&\diiv \bm{u}_b = 0 
&&\text {in } \Omega_b ,\\
&\bm{u}_b= \Kb (- \nabla p_b + \rho_0 \alpha \omega_b \bm{g})
&&\text {in } \Omega_b ,\\
&\nabla_{\tau} \cdot \bm{u_\gamma} = [\lambda] 
&&\text {in } \Gamma^0 ,\\
&\bm{u_\gamma} = b\Kf (-\nabla_{\tau} p_\gamma + \rho_0 \alpha \omega_\gamma \bm{g_\tau}) 
&&\text {in } \Gamma^0 ,\\
& \tr \bm{u}_b \cdot \bm{n}_b = \lambda &&\text{on } \Gamma^\pm, \\
&\lambda = \Kf (- \frac{p_\gamma - \tr p_b}{b/2} + \rho_0 \alpha \tr \omega_b g_n) \qquad
&&\text {on } \Gamma^\pm.\\
\end{aligned}
\right.
\end{equation}

In writing (\ref{eq:mixed_darcy}), we have introduced the new variable $\lambda \, \mathrm{[L/T]} $ and the jump operator $ [v] = v^+ + v^- $. 
For $v$ defined on $\Gamma^\pm$, we will consider $[v]$ to be defined on $\Gamma^0$.
Note that $\Sigma$ collapsed onto a lower-dimensional object (the boundary of the center plane $\partial \Gamma^0$) onto which we will have to impose boundary conditions. We will ignore its contribution to the mass exchange between fracture and bulk by setting zero normal flux (note that this flux scales linearly with the fracture aperture $b$) on immersed fracture boundaries (or tips) whereas fractures will inherit boundary conditions from the bulk if they touch the boundary.

Note that two sources of modeling error associated with dimensional reduction: (\emph{i}) by collapsing the thin dimension of the fracture we have reassigned part of the domain which previously belonged to the fracture to the bulk medium and (\emph{ii}) the flux $\lambda$ exchanged between the bulk medium and the fracture is a first order approximation to the true flux due to the presence of the normal gradient of the pressure $\nabla p$. 

\subsection{Mixed-dimensional transport}

We can now re-apply the reduction procedure to the transport problem. We start by splitting the equation in the two domains and prescribing compatibility conditions at the interface
\begin{equation}
\left\{
\begin{aligned}
&\partial_t \omega_b + \nabla \cdot \bm{q_b} = 0 
&& \text {in } \Omega_b , \\
&\bm{q_b}= - D \nabla \omega_b + \omega_b \bm{u}_b
&& \text {in } \Omega_b ,  \\[1mm]
&\partial_t \omega_f + \nabla \cdot \bm{q_f} = 0 
&& \text {in } \Omega_f ,  \\
&\bm{q_f}= - D \nabla \omega_f + \omega_f \bm{u}_f  \qquad
&& \text {in } \Omega_f ,  \\[1mm]
&\tr \omega_b = \tr \omega_f 
&& \text {on } \Gamma ,  \\
&\tr \bm{q}_b \cdot \bm{n}_b = \tr \bm{q_f} \cdot \bm{n}_b 
&& \text {on } \Gamma . 
\end{aligned}
\right.
\label{eq:diffusione-comp}
\end{equation}
We integrate the conservation equation over the fracture, splitting the fluxes into their normal and tangential parts:
\begin{align*}
0
= \int_{-b/2}^{b/2} \partial_t \omega_f + \nabla \cdot \bm{q_f}
= b \partial_t \omega_\gamma + \nabla_\tau \cdot \bm{q_\gamma} - (\bm{q}_b \cdot \bm{n}_b)^+ - (\bm{q}_b \cdot \bm{n}_b)^- = 0 .
\end{align*}
As we did for the fluid mass conservation equation, we introduced fracture quantities $\omega_\gamma = \frac{1}{b} \int_{-b/2}^{b/2} \omega_f $ and the integrated tangential solute mass flux $q_\gamma = \int_{-b/2}^{b/2} \mathrm{T} \bm{q}_f$. As done before, scalar quantities are averaged while vector quantities are integrated.

Unlike the continuity equation, a new approximation is needed in the averaging step to deal with the nonlinearity of the advective term: let us consider a splitting of $\omega_f$ into $\omega_\gamma + \tilde{\omega}_f$, where $\omega_\gamma$ is constant across the fracture and $\tilde{\omega}_f$ is a null average fluctuation, and similarly for $\bm{u}_f$. We will neglect the product of fluctuations assuming that they are small, i.e. 
\begin{equation*}
\begin{aligned}
&\int_{-b/2}^{b/2} \diiv (\omega_f \bm{u}_f)
\approx \omega_\gamma \bm{u}_\gamma.
%&\begin{split}
%\int_{-b/2}^{b/2} \omega_f \mathrm{T} \bm{u}_f 
%&= \int_{-b/2}^{b/2} (\omega_\gamma + \tilde{\omega}_f) (\frac{1}{b} \bm{u}_\gamma + \mathrm{T} \tilde{\bm{u}}_f) \\
%&= \int_{-b/2}^{b/2} \omega_\gamma \frac{1}{b} \bm{u}_\gamma
%+ \cancelto{0}{ \int_{-b/2}^{b/2} \omega_\gamma \mathrm{T} \tilde{\bm{u}}_f }
%+ \cancelto{0}{ \int_{-b/2}^{b/2} \tilde{\omega}_f \frac{1}{b} \bm{u}_\gamma }
%+ \int_{-b/2}^{b/2} \tilde{\omega}_f \mathrm{T} \tilde{\bm{u}}_f
%\approx \omega_\gamma \bm{u}_\gamma .
%\end{split}
\end{aligned}
\end{equation*}

For fractures of non-negligible aperture however, internal motions (see e.g. intrafractrure mode 2A in figure \ref{fig:3d-modes}) could make this approximation inappropriate.

Just as for Darcy's law, Fick's law will look identical when projected in the fracture plane, while for the normal projection, we will have to resort to a first order approximation
$$
(\tr \bm{q}_b \cdot \bm{n}_b)^\pm = D (\frac{\omega_b^\pm - \omega_\gamma}{b/2}) + \omega_b^\pm \lambda^\pm
$$

\begin{figure}
\centering
\includegraphics[width=0.6\textwidth]{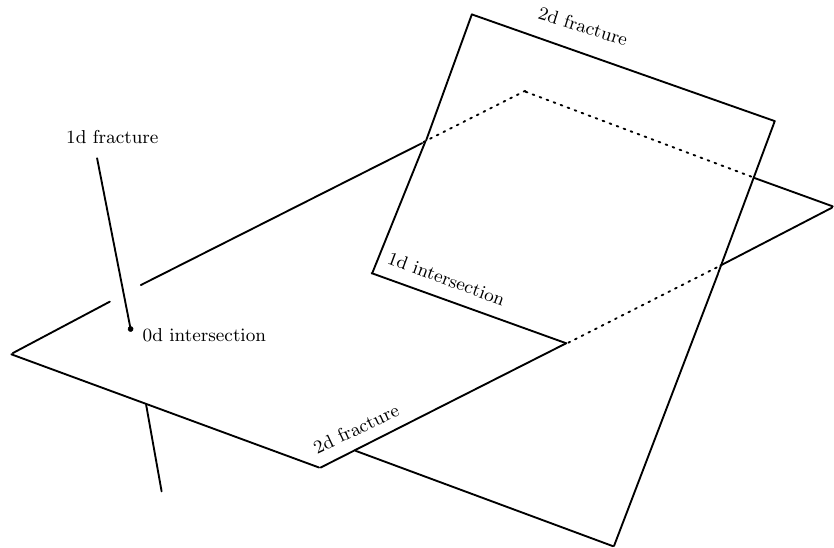}
\caption{Multidimensional coupling for a sample fracture network embedded in a three-dimensional domain.}
\label{fig:intersezioni}
\end{figure}

Finally, we obtain a mixed-dimensional system for the transport of a solute in a fractured porous medium,

\begin{equation}
\label{eq:single-fracture-transport}
\left\{
\begin{aligned}
&\partial_t \omega_b + \nabla \cdot \bm{q}_b = 0                                                                            
&& \text {in } \Omega_b ,\\
&\bm{q}_b = -D \nabla \omega_b + \omega_b \bm{u}_b                                                                           
&& \text {in } \Omega_b ,\\
&b\, \partial_t \omega_\gamma + \nabla_\tau \cdot \bm{q}_\gamma = [\theta]
&& \text {in } \Gamma^0 ,\\
&\bm{q}_\gamma = -b D \nabla_\tau \omega_\gamma + \omega_\gamma \bm{u}_\gamma
&& \text {in } \Gamma^0 ,\\
& \tr \bm{q}_b \cdot \bm{n}_b = \theta 
&&\text{on } \Gamma^\pm, \\[2mm]
&\theta = D (\frac{\tr \omega_b - \omega_\gamma}{b/2}) + \tr \omega_b \lambda                                                
&& \text {on } \Gamma^\pm ,\\
\end{aligned}
\right.
\end{equation}

to be complemented with boundary and initial conditions.
It is important to notice that the equations in the bulk and the equations in the fracture, both in \eqref{eq:mixed_darcy} and \eqref{eq:single-fracture-transport}, are entirely decoupled apart from their interaction through the interface variables $\lambda$ and $\theta$. 
The way these variables appear in the equations reveals how the domains are coupled across different dimensions: in the higher-dimensional domain interface variables appear as Neumann boundary conditions, while in the lower-dimensional fracture they appear as sources inside the domain.
Although the averaging procedure has been carried out for a three-dimensional domain with a single fracture,
the exact same procedure works for a generic $n$-dimensional domain ($ n = 1, 2, 3 $) with multiple, possibly
intersecting fractures. In this case, the procedure is conceptually carried out hierarchically: $n$-dimensional
quantities are coupled to $(n-1)$-dimensional quantities through $(n-1)$-dimensional interface fluxes (figure \ref{fig:intersezioni}).

% SIDE CAPTION FIGURA GRANDE

In view of dealing with the general case of multiple fractures of different dimensions, we want to write equations which hold for each dimension.
Apart from making the system of equations more compact, it will also make the mechanism by which different domains interact more clear.

We collect all domains of equal dimension $d \in \{0, 1, 2, 3\}$ in a single domain $\Omega_d$, and denote by $ \Gamma_d $ the interface between domains $\Omega_d$ and $\Omega_{d+1}$. We also introduce mixed-dimensional variables $ \omega_d, p_d $, fluxes $ \bm{u}_d, \bm{i}_d $ and in-plane gradient operator $\nabla_d$ defined on $\Omega_d$, interface fluxes $ \lambda_d, \theta_d $ defined on $ \Gamma_d $.

With these quantities at hand, we can generalize systems (\ref{eq:mixed_darcy}, \ref{eq:single-fracture-transport}) to 
\begin{equation}
\label{eq:cont_frac}
\left\{
\begin{aligned}
& \nabla_d \cdot \bm{u}_d = [\lambda_d] 
&& \text{in } \Omega_d ,\\
& b^{3-d} \partial_t \omega_d + \nabla_d \cdot \bm{q}_d = [\theta_d]
&& \text{in } \Omega_d ,\\
& \bm{u}_d = b^{3-d} k_d/\phi\mu\; \left(-\nabla p_d + \rho_0 \alpha \omega_d \bm{g}_d^\tau  \right) 
&& \text{in } \Omega_d ,\\ 
& \bm{q}_d = - b^{3-d} D \nabla \omega_d + \omega_d \bm{u}_d 
&& \text{in } \Omega_d ,\\[2mm]
& \tr \bm{u}_d \cdot \bm{n}_d = \lambda_{d-1}
&&\text{on } \Gamma^{d-1}, \\
& \tr \bm{q}_b \cdot \bm{n}_b = \theta_{d-1}
&&\text{on } \Gamma^{d-1}, \\[2mm]
& \lambda_d = b^{2-d} k_d/\phi\mu\; \left(\frac{ \tr p_{d+1} - p_d }{b/2} + \rho_0 \alpha \omega_d \bm{g}_d^n \right)
\qquad\qquad && \text{on } \Gamma_d ,\\
& \theta_d = b^{2-d} D \frac{ \tr \omega_{d+1} - \omega_d }{b/2} + \omega_d \lambda_d 
&& \text{on } \Gamma_d .\\
\end{aligned}
\right.
\end{equation}

with suitable boundary conditions.

\subsection{The impact of fractures on convection onset}
\label{sec:frac-hrl}

The original discussion of the HRL problem addressed the question of the possible onset of convection. 
This was expressed as a critical Rayleigh number $\Ra_c$, such that convective motion is possible for $\Ra \geq \Ra_c$.

In recent years, different studies have tried to understand in what way fractures influence the possibility of convection.
It is clear that highly permeable fractures constitute preferential paths for flow, thus enabling convection or enhancing its strength.
As shown in \cite{vg14}, for large fracture density, calculating an average Rayleigh number based on the average (upscaled) permeability (neglecting the specific fracture configuration) can be effective at predicting the onset and strength of convection. 

\begin{figure}
\centering\begin{minipage}[t]{0.5\textwidth}\hfill\raisebox{-\height}{
\includegraphics{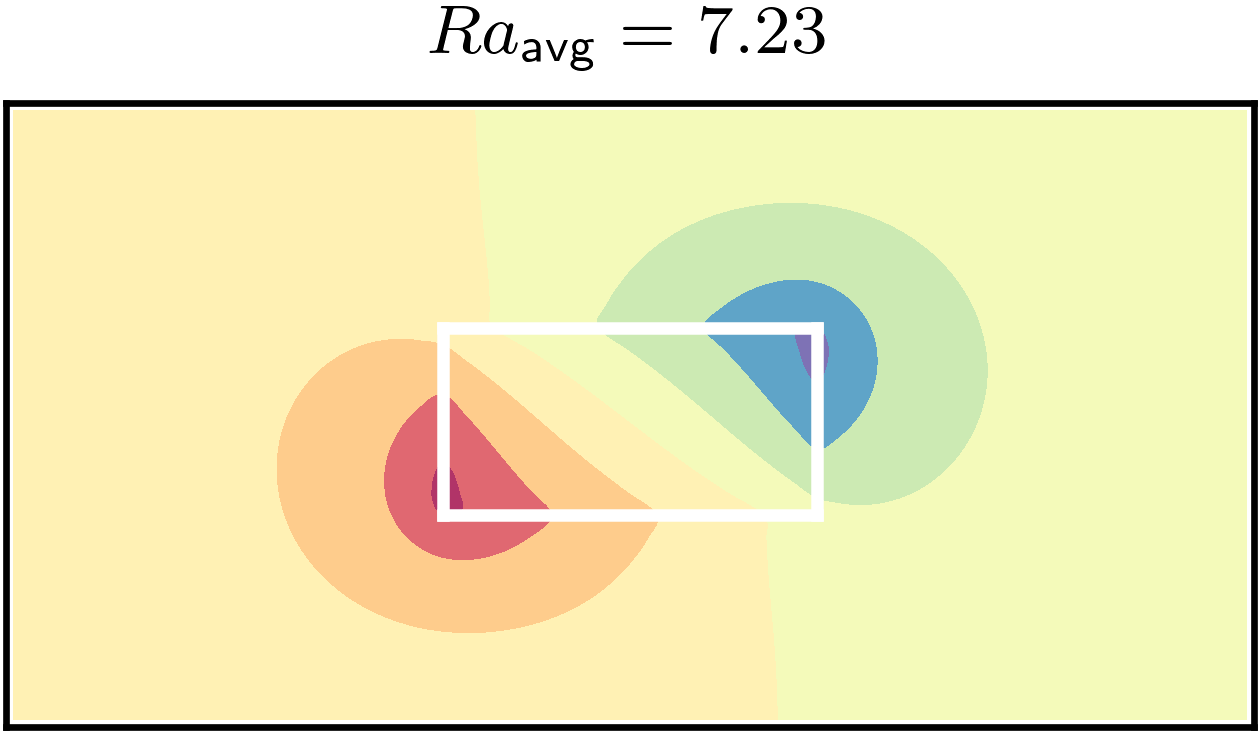}
}\end{minipage}\hspace{0.04\textwidth}\begin{minipage}[t]{0.4\textwidth}
\caption{In this configuration the Rayleigh number computed based on the upscaled permeability is well below the critical Rayleigh number. The presence of the fracture circuit however enables convection regardless.}
\label{fig:bulk-cond}
\end{minipage}
\end{figure}

For lower fracture densities however, this approach is not adequate.
As shown in figure \ref{fig:bulk-cond}, a continuous fracture loop barely modifying the permeability, such that the average Rayleigh number is well below the critical Rayleigh number for homogeneous media, still exhibits convective motion.
In \cite{vg14}, the key factor for enabling convection in the case of low-density fracture configurations is shown to be the presence of continuous fracture circuits, around which convection cells can form.
Simple scenarios (also used as validation cases in this work) tried to relate the location, aspect ratio and size of fracture circuits to the possibility and strength of convection.

A more systematic study remains to be done, with the aim of uncovering better quantitative relations that can be extended to more complex fracture configurations. Moreover, we have to consider that in applicative scenarios the particular geometry of the underground fracture network is mostly unknown, or in the best cases described by statistical parameters such as fracture density, mean lengths and orientations.
In these cases, being able to relate statistical parameters such as the ones mentioned to a quantitative estimate on the strength of convection and its uncertainty would be both useful from an applicative standpoint and interesting in its own right. 

% SIDE CAPTION FIGURA PICCOLA

\begin{remark}
Let us consider the diffusion of a solute in a domain cut by a horizontal fracture.
The continuous equidimensional problem (\ref{eq:diffusione-comp}) admits the linear concentration profile solution for the boundary conditions prescribed by the HRL problem.
$$
\omega_j = \omega_{max} \frac{y}{H},\quad \bm{i}_j = -D \nabla \omega_j = -D \frac{\omega_{max}}{H} \bm{e_y}
$$
for $ j \in \{ B, f \} $.

With dimensional reduction, the system of equations to be solved is replaced by (\ref{eq:single-fracture-transport}), which yields a piecewise linear concentration profile:
\begin{align*}
&\omega_b = 
    \left\{
    \begin{aligned}
    &(\omega_{max} - \delta \omega) \, {y}/{H} \qquad && y < H/2 , \\
    &(\omega_{max} - \delta \omega) \, {y}/{H} + \delta \omega \quad && y > H/2 ,
    \end{aligned}
    \right.
\\
&\omega_f = \omega_{max}/2 \\
&\delta \omega = \frac{\omega_{max}}{1 + H/b} ,
\end{align*}
$b$ being the fracture aperture. Fractures which qualify as thin enough to be treated as lower-dimensional regions will always satisfy $ b \ll H $, 
thus making the error in the concentration profile small: $ \delta \omega \ll \omega_{max} $. This small model reduction error however can manifest itself in the form of (small) artificial fluxes around fracture tips: the concentration gradient that arises from matching the two solutions creates a circulating diffusive flux, as illustrated in figure \ref{fig:flussi-art} on the left. The exact same reasoning can then be applied for the fluid mass conservation equation: small artificial fluid mass fluxes may appear around fracture tips due to the discontinuity of pressure across the fracture (see figure \ref{fig:flussi-art} on the right). 

\begin{figure}
\centering
\includegraphics[width=\textwidth]{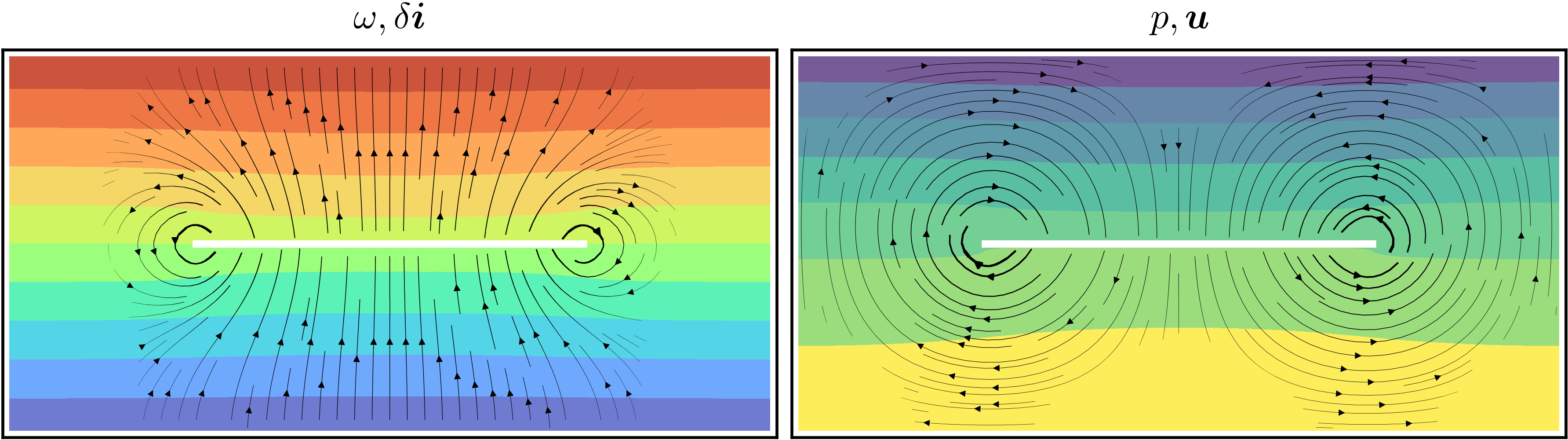}
\caption{ 
Artificial circulating fluxes arising from reducing the fracture to lower-dimensional domains.
Diffusive fluxes on the left are due to concentration discontinuities.
Fluid flow on the right is due to pressure discontinuity.
}
\label{fig:flussi-art}
\end{figure}
\end{remark}

\section{Numerical discretization}
\label{sec:discretization}

This section is dedicated to the discretization of system \eqref{eq:cont_frac}.
The method chosen for the spatial discretization of the system is the finite volume method.
In the fractured problem, we have a sequence of domains $\Omega_d,\, d \in \{0,1,2,3\}$.
We start by defining a mesh on each of them.
While the different meshes in principle can be completely independent, the mathematical formulation is more easily expressed using conforming meshes: lower dimensional meshes are implicitly defined by their higher dimensional neighbour.

\begin{figure}
\centering
\includegraphics[width=\textwidth]{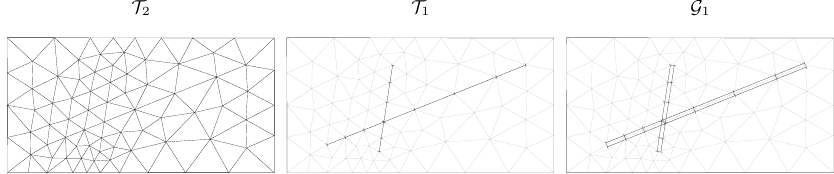}
\caption{Meshes used for correctly formulating the discrete problem. Note how the higher-dimensional mesh conforms to the fracture geometry.}
\end{figure}

The sequence of meshes for each of the domains $\Omega_d$ will be denoted by $\mathcal{T}_d$, the set of edges of $\mathcal{T}_d$ by $\mathcal{E}_d$.
Even though the interfaces $\Gamma_d$ are co-located with the lower dimensional domains $\Omega_d$, we keep the two separate by defining interface meshes $\mathcal{G}_d$.
In what follows, $K$ will denote a generic element of mesh $\mathcal{T}_d$, $\gamma$ a generic element of the interface mesh $\mathcal{G}_d$ and $\sigma$ a generic face of $\mathcal{E}_d$.
We will also indicate the normal pointing inside the fracture at element $\gamma$ as $\bm{n}_\gamma$.

The set $\mathcal{E}_d$, or equivalently the faces of any element $K$, can be partitioned into three sets:
(\emph{i}) internal or belonging to the Dirichlet boundary,
(\emph{ii}) belonging to the Neumann boundary,
(\emph{iii}) adjacent to a lower-dimensional domain.
Notice that we get two different partitions based on which boundary conditions we use to split the boundary (boundary conditions related to the flow problem or to the transport problem): 
$\partial K = \partial K_i^p \cup \partial K_N^p \cup \partial K_f = \partial K_i^\omega \cup \partial K_N^\omega \cup \partial K_f $.
Also, thanks to the conforming mesh hypothesis, faces in $ \partial K_f $ can be identified with elements of the interface mesh $ \mathcal{G}_d $, thus enabling us to legitimately write integrals such as $ \int_\gamma \lambda_d $ where $ \gamma \in \partial K_f $.

\begin{figure}
\centering
\includegraphics[width=\textwidth]{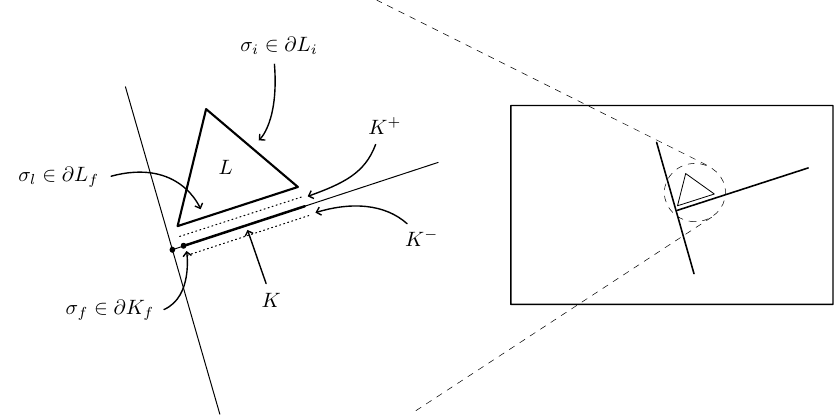}
\caption{Illustration of notation introduced to discretize the mixed-dimensional problem.}
\end{figure}

With the notation in place, we can start integrating the conservation equations in system (\ref{eq:cont_frac}) over a generic element $K \in \mathcal{T}_d$:
\begin{align*}
& \int_K \nabla_d \cdot \bm{u}_d \dd x^d = \int_K [\lambda] \dd x^d , \\
& \int_K b^{3-d} \frac{\partial \omega_d}{\partial t} \dd x^d + \int_K \nabla_d \cdot \bm{q}_d \dd x^d = \int_K [\theta] \dd x^d .
\end{align*}
Working on the integrals one by one we have
\begin{alignat*}{1}
&\int_K \nabla_d \cdot \bm{u}_d \dd x^d 
= \int_{\partial K} \bm{u}_d \cdot \bm{n} \dd x^{d-1}
= \sum_{\partial K_i} \int_{\sigma} \bm{u}_d \cdot \bm{n} \dd x^{d-1}
+ \sum_{\partial K_f} \int_{\sigma} \lambda_{d-1} \dd x^{d-1}
+ \sum_{\partial K_N^p} \int_{\sigma} \bm{u}_N \dd x^{d-1} , \\
& \int_K [\lambda_d] \dd x^d = \int_{K^\pm} \lambda_d \dd x^d , \\
& \int_K b^{3-d} \de_t \omega_d \dd x^d = b^{3-d} \frac{d}{dt} \int_K \omega_d \dd x^d , \\
&\int_K \nabla_d \cdot \bm{q}_d \dd x^d 
= \int_{\partial K} \bm{q}_d \cdot \bm{n} \dd x^{d-1}
= \sum_{\partial K_i} \int_{\sigma} \bm{q}_d \cdot \bm{n} \dd x^{d-1}
+ \sum_{\partial K_f} \int_{\sigma} \theta_{d-1} \dd x^{d-1}
+ \sum_{\partial K_N^\omega} \int_{\sigma} \bm{q}_N \dd x^{d-1} , \\
& \int_K [\theta_d] \dd x^d = \int_{K^\pm} \theta_d \dd x^d . \\
\end{alignat*}

Now introducing the discrete variables and fluxes
\begin{equation}
\label{eq:variabili_discrete}
\begin{aligned}
&P_K = \frac{1}{|K|} \int_K p_d \dd x^d ,
\; && U_{K\sigma} \approx \int_\sigma \bm{u}_d \cdot \bm{n}_K \dd x^{d-1} ,
\; && \Lambda_{K^\pm} = \int_{K^\pm} \lambda_d \dd x^d , \\
&W_K = \frac{1}{|K|} \int_K \omega_d \dd x^d ,
&& Q_{K\sigma} \approx \int_\sigma \bm{q}_d \cdot \bm{n}_K \dd x^{d-1} ,
&& \Theta_{K^\pm} = \int_{K^\pm} \theta_d \dd x^d . \\
\end{aligned}
\end{equation}
we can write the discrete version of the conservation equations in (\ref{eq:cont_frac}):
\begin{equation}
\label{eq:conservazione_discr}
\begin{aligned}
&\sum_{\partial K_i} U_{K\sigma}
+ \sum_{\partial K_f} \Lambda_\gamma
+ \sum_{\partial K_N^p} U_{N,\sigma} = \Lambda_{K^\pm} , \\
& b^{3-d} |K| \frac{d}{dt} W_K
+ \sum_{\partial K_i} Q_{K\sigma} 
+ \sum_{\partial K_f} \Theta_\gamma 
+ \sum_{\partial K_N^\omega} Q_{N,\sigma} 
= \Theta_{K^\pm} .
\end{aligned}
\end{equation}
The last step is the discretization of the constitutive laws for the fluxes.
Integrating the constitutive laws gives
\begin{align*}
&\int_{\sigma} \bm{u}_d \cdot \bm{n} \dd x^{d-1} 
= b^{3-d} \Kd \left( 
    \int_\sigma -\frac{\partial p_d}{\partial n} \dd x^{d-1} 
    + \rho_0 \alpha \,\bm{g}\cdot\bm{n}_\sigma \int_\sigma \omega_d \dd x^{d-1}
\right) , \\
&\int_{\gamma} \lambda_d \, \dd x^d 
= b^{2-d} \Kd \int_\gamma 
    \left( \frac{\tr p_{d+1} - p_d}{b/2} + \rho_0 \alpha \tr \omega_{d+1} \,\bm{g}\cdot\bm{n}_\gamma \right) \dd x^d , \\
&\int_{\sigma} \bm{q}_d \cdot \bm{n} \dd x^{d-1} 
= b^{3-d} \left(
    D \int_\sigma -\frac{\partial \omega_d}{\partial n} \dd x^{d-1}
    + \int_\sigma \omega_d \bm{u_d}\cdot\bm{n} \dd x^{d-1}
\right) , \\
&\int_{\gamma} \theta_d \dd x^d 
= \int_\gamma \left( b^{2-d} D \frac{\tr \omega_{d+1} - \omega_d}{b/2}  + \omega_d \lambda_d \right) \dd x^d . \\
\end{align*}
We rewrite each of these laws in terms of the discrete variables defined in (\ref{eq:variabili_discrete}):
\begin{equation}
\label{eq:costitutive_discr}
\begin{aligned}
&U_{K\sigma} =       b^{3-d} |\sigma| \Kd \left(   {\nabla P}_\sigma           + \rho_0 \alpha \,\bm{g}\cdot\bm{n}_\sigma  W_\sigma \right) , \\
&\Lambda_{K^\pm} =   b^{2-d} |K|      \Kd \left(   \frac{P_{K^\pm} - P_K}{b/2} + \rho_0 \alpha \,\bm{g}\cdot\bm{n}_{K^\pm}  W_{K^\pm} \right) , \\
&Q_{K\sigma} =       b^{3-d} |\sigma| D \,\nabla W_\sigma           + W_\sigma U_{K\sigma} , \\
&\Theta_{K^\pm} =    b^{2-d} |K|      D \frac{W_{K^\pm} - W_K}{b/2} + W_{K^\pm} \Lambda_{K^\pm} ,
\end{aligned}
\end{equation}

where the quantities $\nabla \phi_\sigma$ and $\phi_\sigma$, i.e. gradients and face values, depend on the particular finite volume scheme. We choose the Multipoint Flux Approximation scheme described in \cite{aavatsmark02}, which computes the gradient on a face by considering values of all cells sharing a node with the face.
As for all finite volume schemes, fluid mass and solute mass conservation is guaranteed.
In contrast to a two-point scheme (TPFA) however, the MPFA scheme is consistent on general grids.

Note that we make use of a centered scheme for the advective term in the concentration equation. This choice, unlike upwind, is known to produce numerical oscillations for convection dominated flows. 
In all our numerical experiments however, convection is mild enough for our centered scheme to be numerically stable.
We can easily relate the already introduced Rayleigh number to the Peclet number:
$$
\mathit{Pe} = \frac{u h}{D} = \frac{u H}{D} \frac{h}{H} = \Ra \frac{h}{H}
$$
where $h$ is a characteristic cell diameter, $H$ is a characteristic length of the domain.

In all the the following numerical experiments, sufficiently fine grids will be computationally feasible for $\mathit{Pe}$ to be $O(1)$.
The $\mathit{Pe}$ number is nonetheless numerically monitored in all the following simulations.

\subsection{Time discretization and direct solution method}
\label{sec:direct}

We denote as "direct solution method" the integration of the model equations forward in time starting from a zero-solute or equilibrium initial condition, to assess the possible onset of convection.
Once a steady state has been reached, we must verify whether convective motion is present.
We will use the implicit Euler method for advancing in time for its unconditional stability. An adaptive time-stepping is also used to reduce the amount of computation necessary to reach steady state, see section \ref{sec:exp-hrl} for details.

To discretize in time, define a set of timesteps 
$\{t^n\}_{n=0\dotsc N}$
and write our discretized system as:
\begin{equation}
\label{eq:conservazione_discr_t}
\begin{aligned}
&\sum_{\partial K_i} U_{K\sigma}^{n+1}
+ \sum_{\partial K_f} \Lambda_\gamma^{n+1}
+ \sum_{\partial K_N^p} U_{N,\sigma}^{n+1} = \Lambda_{K^\pm}^{n+1} , \\
& b^{3-d} |K| \frac{W_K^{n+1}-W_K^{n}}{\Delta t^{n}}
+ \sum_{\partial K_i} Q_{K\sigma}^{n+1} 
+ \sum_{\partial K_f} \Theta_\gamma^{n+1}
+ \sum_{\partial K_N^\omega} Q_{N,\sigma}^{n+1} 
= \Theta_{K^\pm}^{n+1},\\ 
\end{aligned}
\end{equation}

where $K\in \mathcal{T}_d, \, d\in \{0,1,2,3\}$, $\Delta t^n = t^{n+1}-t^n$, $n\in \{0, \dotsc, N-1\}$.

If the timestep is long enough and the solution stops changing (according to the norm of the difference between two timesteps of the solution), we declare the system to have reached steady state.
Note that since the system of equations is nonlinear, each timestep requires the solution of a nonlinear problem. 
In our case we use Newton iterations by leveraging the automatic differentiation capabilities of the implementation. 
In the test cases, the tolerance for the Newton procedure is set to \num{1e-8} for the concentration increment.

\subsection{Eigenvalue analysis}
\label{sec:eig}

The previously outlined method for assessing stability, while effective at predicting stability, has some shortcomings.
Among the disadvantages are its reliance on the choice of perturbation for the hydrostatic solution (if not starting from zero solute everywhere) and having to reach and assess the steadiness of the solution.
The computational cost of reaching the steady state may become an issue: advancing a nonlinear equation in time with an implicit scheme requires the solution of multiple linear systems for each timestep advancement. 

An alternative method, presented in detail below, relies on inspecting the nonlinear system of equations linearized around the equilibrium solution. 
The nonlinear discrete system (\ref{eq:conservazione_discr}) and (\ref{eq:costitutive_discr}) can be written abstractly as

\begin{equation}
\label{eq:eig_generico}
\begin{aligned}
& M \frac {d W}{d t} + F(W, Y)=0 , \\  
& G(W, Y)=0 ,
\end{aligned}
\end{equation}
where  $W \in \mathbb{R}^n $ collects the degrees of freedom relative to the discrete variable $W$, and $ Y \in \mathbb{R}^{N-n} $ collects the degrees of freedom relative to the discrete variables $P$, $ \Lambda $ and $ \Theta $. 
$F: \mathbb{R}^n \times \mathbb{R}^{N-n} \to \mathbb{R}^n$ and $G: \mathbb{R}^n \times \mathbb{R}^{N-n} \to \mathbb{R}^{N-n}$ collect the linear and nonlinear discrete operators in (\ref{eq:conservazione_discr}).
Any discrete equilibrium solution $(W_s, Y_s)$ will satisfy the system
\begin{align*}
&F(W_s, Y_s)=0 , \\  
&G(W_s, Y_s)=0 .
\end{align*}
To assess  whether the equilibrium solution is also asymptotically stable we perturb the time-dependent system:

\begin{align*}
M (\frac {d W_s}{d t} + \frac {d \delta W}{d t} )+& F(W_s+ \delta W, Y_s + \delta Y)=0 , \\  
& G(W_s+ \delta W, Y_s + \delta Y)=0 ,
\end{align*}
and linearize, taking advantage of the fact that $\delta W, \delta Y$ are small perturbations:
\begin{align*}
M \frac {d W_s}{d t} + M \frac {d \delta W}{d t} +& \cancelto{0}{ F(W_s, Y_s)} + \frac{\partial F}{\partial W} \Bigr|_{W_s, Y_s} \delta W + \frac{\partial F}{\partial Y}\Bigr|_{W_s, Y_s} \delta Y=0 , \\  
& \cancelto{0}{ G(W_s, Y_s)} + \frac{\partial G}{\partial W} \Bigr|_{W_s, Y_s} \delta W + \frac{\partial G}{\partial Y}\Bigr|_{W_s, Y_s} \delta Y=0 .
\end{align*}
Renaming the partial derivatives to ease the notation 
\begin{alignat*}{2}
& A_{ww} = \frac{\partial F}{\partial W} \Bigr|_{W_s, Y_s} , &\qquad& A_{wy} = \frac{\partial F}{\partial Y} \Bigr|_{W_s, Y_s} , \\
& A_{yw} = \frac{\partial G}{\partial W} \Bigr|_{W_s, Y_s} , &&       A_{yy} = \frac{\partial G}{\partial Y} \Bigr|_{W_s, Y_s} , \\
\end{alignat*}
the system becomes
\begin{alignat*}{3}
M \frac {\partial \delta W}{\partial t} + 
&A_{ww} \delta W &&+ A_{wy} \delta Y &&= 0 , \\ 
&A_{yw} \delta W &&+ A_{yy} \delta Y &&= 0 .
\end{alignat*}
Now, relying on the invertibility of $ A_{yy} $ , we can eliminate $\delta Y$ to obtain a single  evolution equation for the perturbation $\delta W$:
\begin{equation*}
M \frac {\partial \delta W}{\partial t} = (A_{wy} A_{yy}^{-1} A_{yw} - A_{ww}) \delta W. \end{equation*}
Substituting an exponential in time trial solution $\delta W(t) = w \, e^{\lambda t} $ in the evolution equation yields an eigenvalue problem:
\begin{equation*}
\lambda M w = (A_{wy} A_{yy}^{-1} A_{yw} - A_{ww}) w,
\end{equation*}
or, by defining the matrix $S= M^{-1} (A_{wy} A_{yy}^{-1} A_{yw} - A_{ww})$,
\begin{equation}
\label{eq:eig}
 S w = \lambda w .
\end{equation}
Note that $\lambda$ describes the evolution of the perturbation in time, while vector $w$ its shape in space because each component represent the value in a grid cell. 
The equilibrium solution $W_s$ is then linearly stable if and only if all the eigenvalues associated to system (\ref{eq:eig}) have negative real part. Conversely, if we can find one or more eigenvalues with positive real part the perturbation $w$ will grow and be sustained in time.

Three computational considerations: (\emph{i}) to assess stability there is no need to compute the entire spectrum, it is enough to compute the eigenvalue of largest real part (\emph{ii}) automatic differentiation of the discrete system (\ref{eq:eig_generico}) can yield the numerical value of the matrices $A_{ww}, A_{wy}, A_{yw}, A_{yy}$ without having to write explicit expressions for them (\emph{iii}) iterative methods for computation of eigenvalue spectra are available e.g. power iterations which do not require explicit expression for the matrix under study, only the ability of computing matrix-vector products. In our case, this removes the need of explicitly inverting the matrix $A_{yy} $. Furthermore, since the matrix only depends on the equilibrium solution, we can factorize the matrix only once using e.g. LU decomposition for fast matrix-vector products during the computation of the spectrum.

\subsection{Implementation details}

The implementation of the numerical methods outlined above is based on the \texttt{PorePy} library \cite{porepy}, which provides the necessary tools for meshing fractured domains and assembling the discrete mixed-dimensional operators. The framework also implements forward automatic differentiation, providing the numerical Jacobians used for performing Newton iterations.

For the eigenvalue analysis (\ref{sec:eig}), the automatic differentiation part provides the Jacobian matrix.
Once the different matrix blocks are identified, we can easily define the matrix-vector product procedure yielding $ S v $.
For the computation of few leading eigenpairs, the Krylov-Schur algorithm (outlined in \cite{stewart02}) has been combined with the dynamic restarting scheme described in \cite{stathopoulos98}. 
Using the inner product defined by the mass matrix for the orthogonalization part of the algorithm has been particularly beneficial in accelerating convergence. 
The use of a custom procedure for computing eigenvalues has been preferred to a default implementation such as \texttt{ARPACK} mostly due to the possibility of monitoring convergence and as a possible starting point for devising more efficient algorithms.

\section{Results}

The model and its numerical approximation have been validated against three reference papers: \cite{diersch02} which treats the Elder problem, \cite{vg14} and \cite{vg15} that treat the HRL scenario respectively in two and three dimensions. Both problems have been extensively used in the literature as benchmarks in the context of density driven flows.

\subsection{Elder problem}
\label{sec:exp-elder}

The Elder problem was originally proposed in a paper by Elder \cite{elder67}, studying thermal convection in a Hele-Shaw cell.
It was later reformulated into a solute convection problem by Voss and Souza \cite{voss1987variable},
where the system of equations is similar to the homogeneous problem \eqref{eq:homogeneous}.

What this benchmark case aims to highlight and validate is the possibility of flow driven purely by density differences: no pressure gradient is being enforced by the boundary conditions.

The problem we want to solve is (\ref{eq:homogeneous}), the domain being \( \Omega = [0, 600] \times [0,150] \). Boundary conditions are
\[
\left\{
\begin{aligned}
    & \bm{u} \cdot \bm{n} = 0 & \quad & \text{on } \partial\Omega ,\\
    & \omega = \omega_{max} && \text{on } \partial\Omega_i ,\\
    & \omega = 0 && \text{on } \partial\Omega_o ,\\
    & \bm{q} \cdot \bm{n} = 0 & \quad & \text{on } \partial\Omega \setminus (\partial\Omega_i \cup \partial\Omega_o) ,\\
\end{aligned}
\right.
\]
where \( \partial\Omega_i = (150, 450) \times 150 \) and \( \partial\Omega_o = (0, 600) \times 0 \).
Since the boundary conditions for the Darcy problem are of Neumann type over the whole boundary we have an ill-posed problem; 
we can however restore the well-posedness by adding an additional constraint, such as imposing zero mean pressure over the whole domain:
\( \int_\Omega p = 0 \).

Initial conditions prescribe \( \omega(x, 0) = 0 ,\, x \in \Omega \).
Equations are integrated in time until $ T = \SI{20}{yr}  $.

\begin{figure}
\centering 
\includegraphics[width=0.95\textwidth]{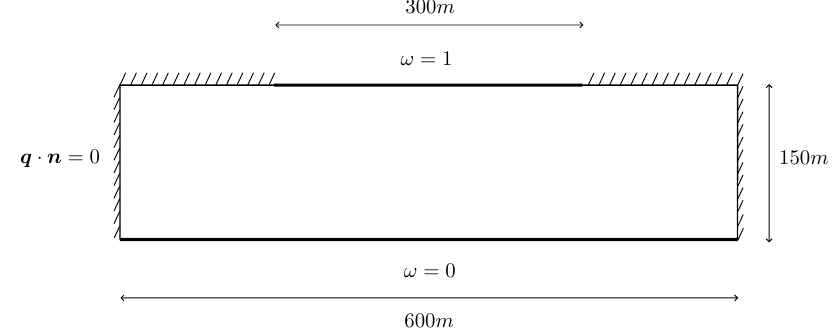}
\caption{Concentration boundary conditions for the Elder problem. Pressure boundary condition are of no flow all around: $ \bm{u}\cdot\bm{n} = 0 $}
% \label{fig:elder}
\end{figure}

All the other parameters of the problem are reported in table \ref{table:elder}.

\begin{table}
\small
\centering 
\begin{tabular}{ l l l l}
% \begin{tabularx}{0.8\linewidth}{ l l l}
\toprule
Permeability                                  & $k$              & \num{4.845e-13}           &  \si{\meter^2} \\
Porosity                                      & $\phi$           & 0.1                       &  1  \\
Viscosity                                     & $\mu$            & \num{1e-3}                &  \si{\kg\per\meter\per\second} \\
Freshwater density                            & $\rho_0$         & \num{1000}                &  \si{\kg\per\meter^3} \\
Solute expansion coefficient                  & $\alpha$         & \num{0.2}                 &  1  \\
Maximum concentration                         & $\omega_{max}$   & 1                         &  1  \\
Gravitational acceleration \qquad\qquad       & $g$ \qquad       & 9.81                      &  \si{\meter\per\second^2} \\
Diffusivity                                   & $D$              & \num{3.565e-6}            &  \si{\meter^2\per\second} \\
\bottomrule 
\end{tabular}
\caption{Parameters of the Elder problem.}
\label{table:elder}
\end{table}

In the absence of gravity, and for small enough $\Ra$, the solute will diffuse from the inlet \( \partial \Omega_i \) until the solution reaches the diffusive steady state. The characteristic time of the evolution is \( T_\text{diff} = H^2/D \), which, for the parameters listed above is 
\( T = (\SI{150}{\meter})^2/\SI{3.565e-6}{\meter^2\per\second} \approx \SI{200}{yr} \).

The corresponding characteristic time associated to convection induced by density differences is much smaller: 
\[T_\text{adv} = \frac{H\phi\mu}{k \rho_0 \alpha \omega_{max} g} \approx \SI{0.5}{yr}. \]
Note that the Rayleigh number presented above as a ratio of velocities \ref{eq:rayleigh} can be equivalently interpreted as a ratio of these timescales:
\[ \Ra = \frac{T_\text{diff}}{T_\text{adv}} \approx 400 .\]

As noted in \cite{diersch02}, the solutions of the Elder problem present grid dependence even for moderately fine meshes, also due to the possible non-uniqueness of the solution.
For this reason, instead of trying to achieve grid independence, we directly validate the implementation by comparing the solutions with \cite{diersch02} for corresponding levels of grid refinement.

\begin{figure}
    \centering
    \includegraphics[width=\textwidth]{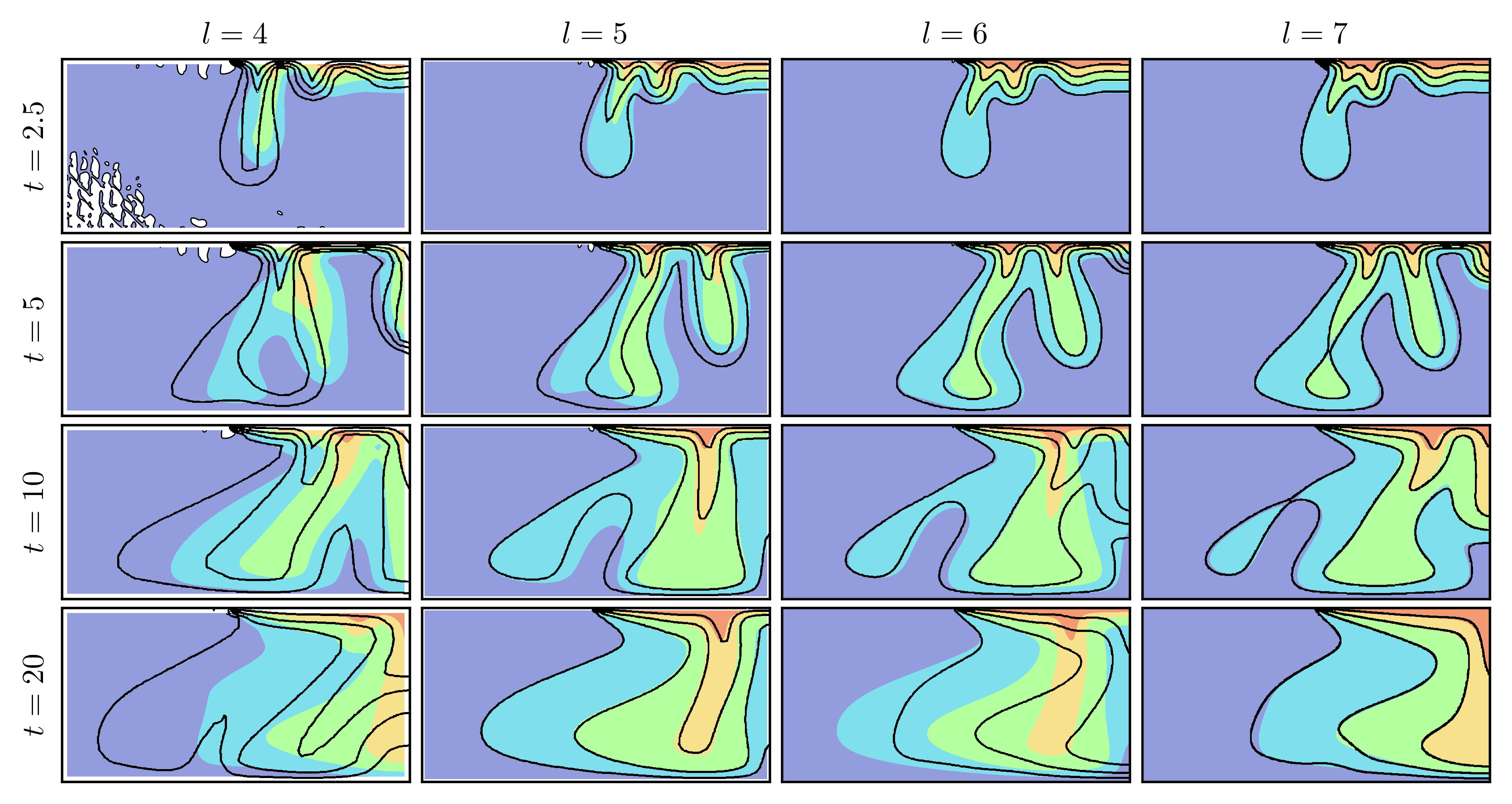}
    \caption{Concentration profile for varying grid refinement at different times.
    Color plot is the result of our implementation.
    Black contour lines are from \cite{diersch02}. 
    White spots indicate negative concentrations (numerical artifact). Timestep is fixed at $ \Delta t = \SI[parse-numbers=false]{1/12}{yr} $, number of (quadrilateral) 
    elements $\mathit{NoE} = 2^{2l + 1} $. }
    \label{fig:elder_val}
\end{figure}

We use quadrilateral grids identical to the ones reported in the reference paper: number of cells $ n =  2^{2l + 1} ,\, l \in \{ 4, 5, 6 \} $ and fixed timestep $ \Delta t = \SI[parse-numbers=false]{1/12}{yr} $. The qualitative comparison using the contours of the concentration profiles is reported in figure \ref{fig:elder_val}.

While the differences in the continuous model (which is not explicitly detailed in \cite{diersch02})
and discretization method (\cite{diersch02} uses adaptive time stepping and Galerkin-FEM) cause different solutions for low levels of grid refinement,
for higher levels of grid refinement the solutions are in good agreement.

\subsection{HRL problem}
\label{sec:exp-hrl}

\begin{figure}
\centering
\includegraphics[width=0.5\textwidth]{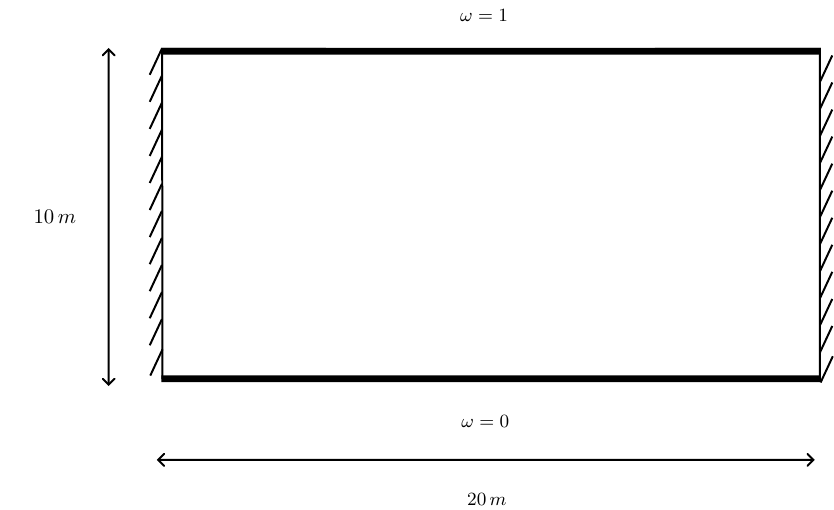}
\caption{Concentration boundary conditions for the HRL problem. Pressure boundary condition are of no flow all around: $ \bm{u}\cdot\bm{n} = 0. $}
\label{fig:vg14}
\end{figure}

The HRL problem, already introduced in sections \ref{sec:hrl} and \ref{sec:frac-hrl}, is the test case analyzed in \cite{vg14}, 
which we will use as reference solutions for validating both the direct method and the eigenvalue method. 
The geometry of the domain and the parameters common to the different simulations are reported in figure \ref{fig:vg14} and table \ref{table:vg14_param}.

\begin{table}
\small
\centering 
\begin{tabular}{ l c c l }
    % \begin{tabularx}{0.8\linewidth}{ l l l}
\toprule
Permeability                                  & $k$              & \num{1e-16}               &  \si{\meter^2} \\
Porosity                                      & $\phi$           & 0.1                       &  1  \\
Viscosity                                     & $\mu$            & \num{1.1e-3}              &  \si{\kg\per\meter\per\second} \\
Freshwater density                            & $\rho_0$         & \num{1000}                &  \si{\kg\per\meter^3} \\
Solute expansion coefficient                  & $\alpha$         & \num{0.7}                 &  1  \\
Maximum concentration                         & $\omega_{max}$   & 0.1                       &  1  \\
Gravitational acceleration \qquad\qquad       & $g$ \qquad       & 9.81                      &  \si{\meter\per\second^2} \\
Diffusivity                                   & $D$              & \num{1e-9}                &  \si{\meter^2\per\second} \\
\bottomrule
\end{tabular}
\caption{Parameters for the HRL problem. Taken from \cite{vg14}. }
\label{table:vg14_param}
\end{table}

The system of equations solved is system (\ref{eq:cont_frac}), with boundary conditions given by
\[
\left\{
\begin{aligned}
    & \bm{u} \cdot \bm{n} = 0 & \quad & \text{on } \partial\Omega ,\\
    & \omega = \omega_{max} && \text{on } \partial\Omega_i ,\\
    & \omega = 0 && \text{on } \partial\Omega_o ,\\
    & \bm{q} \cdot \bm{n} = 0 & \quad & \text{on } \partial\Omega \setminus (\partial\Omega_i \cup \partial\Omega_o) ,\\
\end{aligned}
\right.
\]
where \( \partial\Omega_i = (0, 20) \times 10 \) and \( \partial\Omega_o = (0, 20) \times 0 \).
As for the Elder problem, the model is supplemented by the additional constraint $ \int_\Omega p = 0 $ to obtain a well-posed problem.
The solution strategy is inspired by the one outlined in the reference paper: the initial diffusive steady state concentration is perturbed
and the solution is advanced in time, gradually increasing the timestep $ \Delta t $.
Criteria for adapting the timestep $\Delta t$ include both the number of Newton iterations required for convergence and 
the norm of the concentration difference $ \Vert \omega^{n+1} - \omega^{n} \Vert $.
If $ \Delta t $ is large enough with respect to characteristic time scales, and the concentration difference between timesteps is small enough, we consider the solution to have reached steady state.

\begin{figure}
    \centering
    \includegraphics{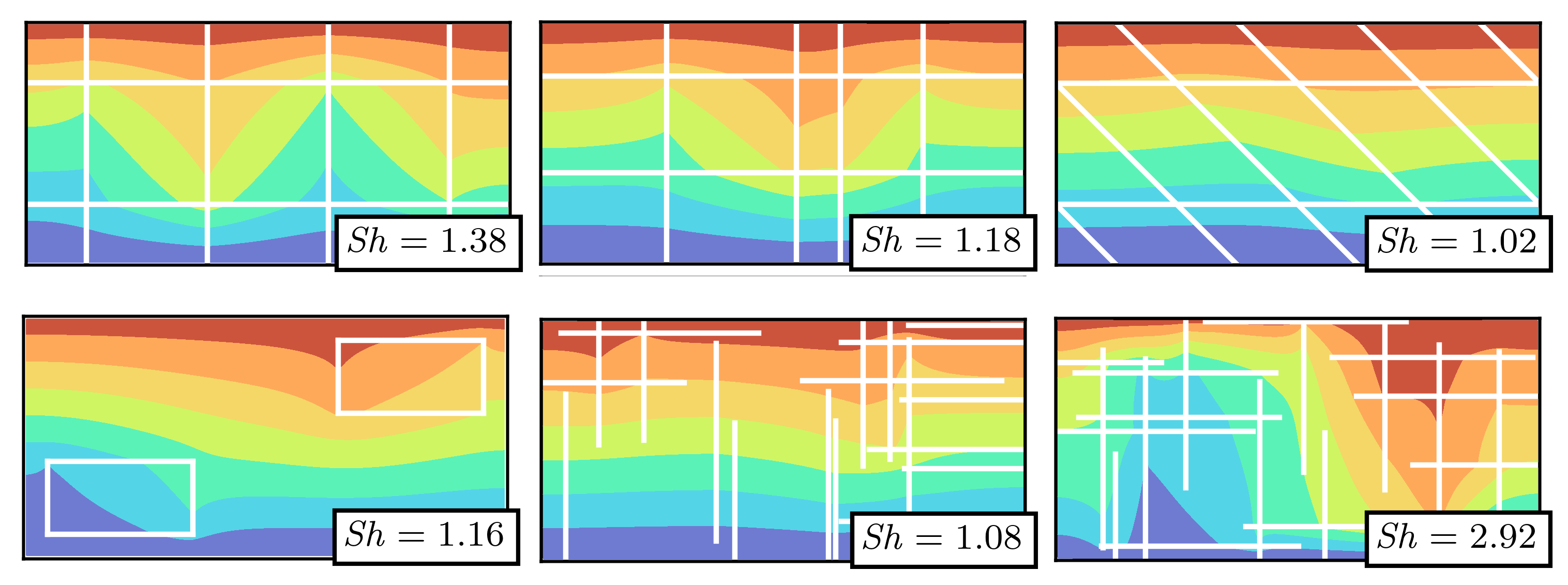}
    \caption{Concentration profiles for a few selected test cases from \cite{vg14}.} 
    \label{fig:validation_vg14}
\end{figure}

A quantitative comparison between our results and the one presented in the reference paper is presented in table \ref{table:validation_vg14t}.
Concentration profiles for a few selected cases are also presented in figure \ref{fig:validation_vg14}.
In all the test cases, there is agreement between our results and \cite{vg14} on whether convective motion is present at steady state.
In the majority of cases, there is both qualitative agreement in the concentration profiles and quantitative agreement on the strength of convection, as measured by the Sherwood number, even in test cases with complex fracture configurations such as $\mathsf{E9a}$ and $\mathsf{E9b}$. 
As for the significant differences in cases such as $ \mathsf{B3}, \mathsf{C2}, \mathsf{C3} $
they might be due to important differences in the model solved in this study and the one in the reference paper:
we use the Boussinesq approximation, neglect dispersivity, and use a different numerical method. 
These differences however seem to have a modest overall impact.

\begin{table}
\small\centering 
\begin{tabular}{ cc }
\begin{tabular}{ l r r r }\toprule 
Scenario & $ \Sh_\text{vg} $ & $ \Sh $ & $ \epsilon_S $ \\[1pt]
\midrule
$\mathsf{A1}$    &  1.00 &  1.00 &     0.00  \\
$\mathsf{A2}$    &  1.36 &  1.38 &     0.01  \\
$\mathsf{A3}$    &  1.56 &  1.77 &     0.14  \\
$\mathsf{A4}$    &  1.75 &  1.79 &     0.02  \\
$\mathsf{B1}$    &  1.00 &  1.00 &     0.00  \\
$\mathsf{B2}$    &  1.13 &  1.18 &     0.04  \\
$\mathsf{B3}$    &  1.49 &  1.19 &    -0.20  \\
$\mathsf{B4}$    &  1.32 &  1.29 &    -0.03  \\
$\mathsf{C1}$    &  1.00 &  1.00 &     0.00  \\
$\mathsf{C2}$    &  1.17 &  1.02 &    -0.13  \\
$\mathsf{C3}$    &  1.21 &  1.06 &    -0.13  \\
$\mathsf{C4}$    &  1.21 &  1.13 &    -0.06  \\
$\mathsf{E9a}$   &  1.06 &  1.08 &     0.02  \\
\bottomrule\end{tabular}
\qquad
\begin{tabular}{ l r r r }\toprule
Scenario & $ \Sh_\text{vg} $ & $ \Sh $ & $ \epsilon_S $ \\
\midrule
$\mathsf{D1}$    &  1.08 &  1.08 &     0.00  \\
$\mathsf{D2}$    &  1.00 &  1.00 &     0.00  \\
$\mathsf{D3}$    &  1.00 &  1.00 &     0.00  \\
$\mathsf{D4}$    &  1.01 &  1.02 &     0.01  \\
$\mathsf{D5}$    &  1.07 &  1.06 &    -0.01  \\
$\mathsf{D6}$    &  1.08 &  1.09 &     0.01  \\
$\mathsf{D7}$    &  1.28 &  1.27 &    -0.01  \\
$\mathsf{D8}$    &  1.15 &  1.16 &     0.01  \\
$\mathsf{D9}$    &  1.15 &  1.16 &     0.01  \\
$\mathsf{D10}$   &  1.45 &  1.45 &     0.00  \\
$\mathsf{D11}$   &  1.37 &  1.37 &     0.00  \\
$\mathsf{D12}$   &  1.39 &  1.39 &     0.00  \\
$\mathsf{E9b}$   &  2.69 &  2.92 &     0.09  \\
\bottomrule\end{tabular}
\end{tabular}
\caption{
Quantitative comparison of strength of convection at steady state as measured by the Sherwood number (defined in (\ref{eq:sherwood})). 
$ \epsilon_S = \frac{\Sh - \Sh_\text{vg}}{\Sh_\text{vg}} $ indicates the relative error. Fracture configurations of the different scenarios are taken from \cite{vg14}. }
\label{table:validation_vg14t}
\end{table}

\begin{figure}
    \centering
    \includegraphics[width=\textwidth]{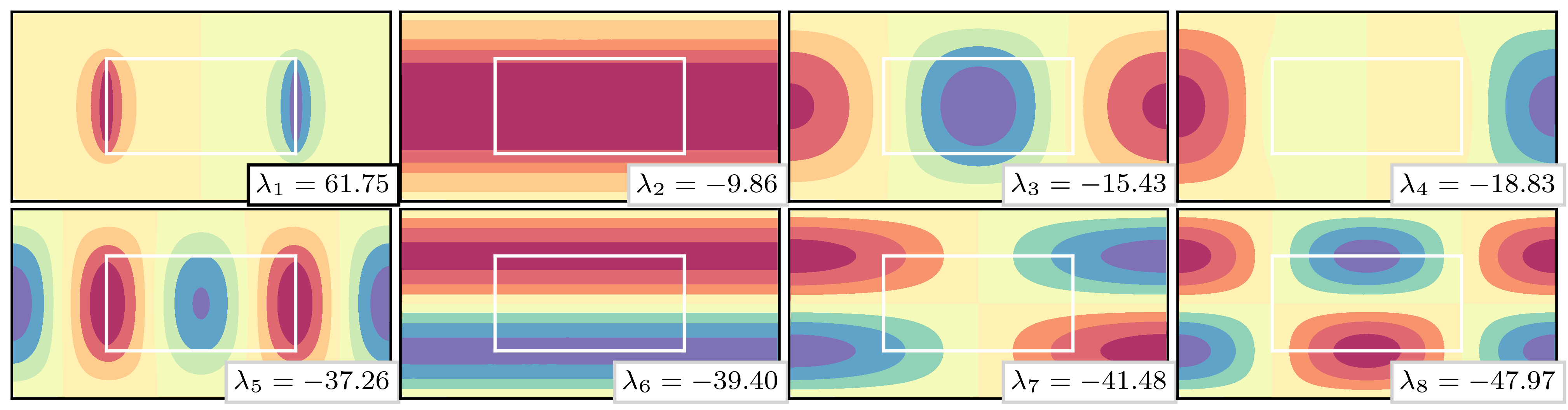}
    \caption{The results of the eigenvalue method applied to scenario $\mathsf{D11}$. The eigenvectors represent perturbations from the equilibrium solution, they thus have no intrinsic scale. Note thatthe first one corresponds to an unstable convection mode.}
    \label{fig:eig_vg14_D11}
\end{figure}

The analysis of the different scenarios presented above can be complemented with the stability analysis based on eigenvalues outlined in section \ref{sec:eig}.
We begin with a detailed analysis of scenario $\mathsf{D11}$.
The method can provide the eigenpairs corresponding to $k$ eigenvalues with largest real part, for reasonably small $k$.
Figure \ref{fig:eig_vg14_D11} shows some of the computed eigenpairs, and, as we can see, one of them is positive, indicating the presence of natural convection. 

The computed eigenvalues are consistent with the direct simulation: indeed, the instability of the diffusive steady state is confirmed by the presence of one eigenvalue of positive real part.
The numerical error due to the iterative nature of the eigenvalue computation is estimated by the formula $ \epsilon = \frac{\lVert S x - \lambda x \lVert}{\lambda x} $, where $S$ is introduced in (\ref{eq:eig}). Grid independence instead is assessed by computing the eigenvalues on a coarser grid and again computing a relative error: $ \epsilon_g = \frac{\left| \lambda - \lambda_g \right|}{\lambda} $, where $ \lambda_g $ are the eigenvalues computed on a coarser grid. 
The corresponding errors for this scenario are reported in table \ref{table:errori_eig_D11}. The errors $\epsilon$, $\epsilon_g$ are overall showing good accuracy for almost every scenario.

\begin{table}
\begin{minipage}{0.5\textwidth}
\small\centering
\begin{tabular}{ r r r }\toprule
\hd{$\lambda$} & \hd{$\epsilon$} & \hd{$\epsilon_g$} \\
% \makecell[c]{ \boldmath\bfseries $\lambda$ } & 
% % \makecell[c]{ \boldmath $\lambda$ } & 
% \makecell[c]{ $\bm{\epsilon}$ }& 
% \makecell[c]{ $\bm{\epsilon_g}$ }\\
\midrule
61.75   &   0.0000 & 0.0014 \\
-9.86   &   0.0001 & 0.0004 \\
-15.43  &   0.0001 & 0.0007 \\
-18.83  &   0.0001 & 0.0007 \\
-37.26  &   0.0001 & 0.0014 \\
-39.40  &   0.0004 & 0.0015 \\
-41.48  &   0.0003 & 0.0016 \\
-47.97  &   0.0008 & 0.0020 \\
\bottomrule
\end{tabular}
\end{minipage}
\begin{minipage}{0.5\textwidth}
\caption{Numerical errors for computed eigenvalues in scenario $\mathsf{D11}$. 
Degrees of freedom ratio for computing grid independence is $ \noe/\noe_g = 1.77 $, where $\noe$ refers to number of mesh elements in the domain (counting both porous medium and fractures). This number is also equal to the size of the eigenvalue problem. } 
\label{table:errori_eig_D11}
\end{minipage}  
\end{table}

Apart from predicting the possibility of convection, the eigenfunctions and eigenvalues can give insight into the evolution dynamics in the vicinity of the diffusive initial condition.
We start by noticing that not only all of the computed eigenvalues have imaginary part equal to zero, but the related eigenfunctions furthermore are all mutually approximately orthogonal. Let $\langle e_i, e_j\rangle$ denote the scalar product between (normalized) eigenfunctions: we have that, if $i\neq j$ the value is about three orders of magnitude smaller than 1 for test case D11. We remark that the $S$ matrix is not symmetric for reasons due to the implementation of boundary conditions, so there is no obvious reason to expect these results.

The orthogonality in particular suggests the possibility of studying the time evolution of the concentration in the space spanned by these eigenfunctions: given the solution $ \omega(t) $ and the eigenfunctions $ \{ e_1, \ldots, e_k \} $, let us define the scalar functions 
\begin{equation}
\alpha_i(t) = \langle \omega(t) - \omega_0, e_i \rangle,\,i=1,\ldots,k ,
\label{eq:alpha}
\end{equation}

representing the projections of $\omega(t)-\omega_0$ on the eigenfunction basis.
At steady state, different eigenfunctions give non-negligible contribution to the solution $\omega(t)$: no obvious relationship exists between the eigenpairs and the steady state solution. This should come as no surprise given that the eigenvalue analysis is localized around the initial condition.
In the early stages of the simulation however, not only the eigenvalue analysis correctly predicts the shape of the growing perturbation (identical to the only eigenfunction associated to a positive eigenvalue), but the eigenvalue $\lambda_1$ also provides a good estimate of its growth rate.

\begin{figure}
\centering
  \includegraphics[width=0.5\textwidth]{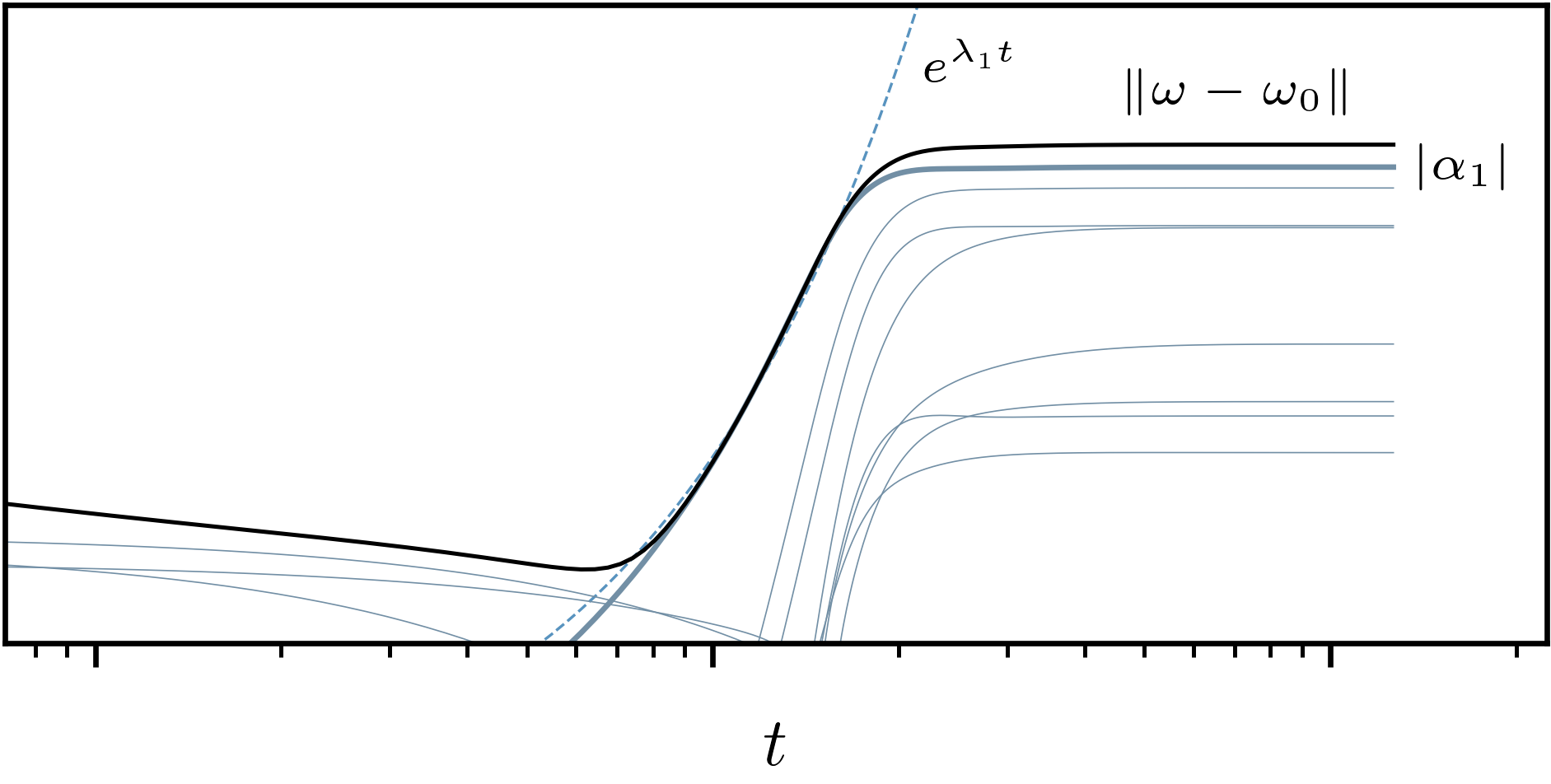}   
\caption{ The norm of the perturbation $ \Vert \omega(t) - \omega_0 \Vert $ and the magnitude of the projections $ \lvert \alpha_j(t) \rvert $ are plotted against time. 
The overlap between $ \Vert \omega(t) - \omega_0 \Vert$ and $\lvert \alpha_1(t) \rvert $ at early times indicates that $ \omega(t) - \omega_0 \propto e_1 $ approximately i.e. the perturbation is growing with the shape of the first eigenfunction.
The overlap between the exponential (dashed line) and the norm of the perturbation indicates that the eigenvalue $\lambda_1 $ is providing a good estimate of the growth rate.}
\end{figure}

Let us now consider other scenarios.  Table \ref{table:eig-vg14} provides the five eigenvalues with largest real part for each of the scenarios from \cite{vg14} presented above. With no exception, a Sherwood number greater than 1 is associated to an eigenvalue with positive real part, thus correctly predicting the possibility or impossibility of convection. Moreover, larger $\Sh$ are (weakly) associated to greater number of eigenvalues or eigenvalues greater in magnitude.

\begin{table}
\small
    \centering
    \begin{tabular}{ c r @{\hskip 0.2in} r r r r r @{\hskip 0.2in} r r c }
    \toprule
    \hd{Scenario} & \hd{$\Sh$} & \hd{$\lambda_1$} & \hd{$\lambda_2$} & \hd{$\lambda_3$} & \hd{$\lambda_4$} & \hd{$\lambda_5$} & \hd{$\epsilon$} & \hd{$\epsilon_g$} & \hd{$\noe/\noe_g$} \\
    \midrule
% 0A    &  1.00 &   -9.87 &  -10.36 &  -14.80 &  -25.24 &  -39.45 &  0.00861 &  0.00058 &  1.78 \\
% 0B    &  1.00 &   -0.21 &   -4.52 &   -5.02 &   -9.87 &  -18.07 &  0.00481 &  0.02957 &  1.78 \\
% 0C    &  1.00 &    0.19 &   -4.36 &   -4.48 &   -9.87 &  -17.44 &  0.00353 &  0.03303 &  1.78 \\
% 0D    &  1.00 &    9.86 &    8.90 &   -0.50 &   -1.97 &   -9.87 &  0.00631 &  0.00563 &  1.78 \\
$\mathsf{A1}$    &  1.00 &   -9.87 &  -10.88 &  -16.61 &  -25.55 &  -39.35 &  0.00919 &  0.00075 &  1.84 \\
$\mathsf{A2}$    &  1.36 &   51.92 &   30.36 &    2.30 &   -9.87 &  -39.43 &  0.01402 &  0.01351 &  1.86 \\
$\mathsf{A3}$    &  1.56 &   16.03 &   13.45 &    8.44 &    0.17 &   -6.44 &  0.03232 &  0.07508 &  1.93 \\
$\mathsf{A4}$    &  1.75 &   11.08 &   10.99 &    0.69 &   -0.17 &   -9.87 &  0.00652 &  0.03490 &  1.47 \\
$\mathsf{B1}$    &  1.00 &   -9.87 &  -10.76 &  -15.70 &  -24.87 &  -39.43 &  0.00500 &  0.00110 &  1.75 \\
$\mathsf{B2}$    &  1.13 &   31.29 &   26.31 &    0.96 &   -9.87 &  -18.21 &  0.00265 &  0.04576 &  1.73 \\
$\mathsf{B3}$    &  1.49 &   15.31 &    9.31 &    1.63 &   -4.05 &   -9.87 &  0.00658 &  0.04265 &  1.64 \\
$\mathsf{B4}$    &  1.32 &   11.62 &    4.39 &   -0.54 &   -3.63 &   -9.87 &  0.02193 &  0.05261 &  1.39 \\
$\mathsf{C1}$    &  1.00 &   -9.87 &  -10.59 &  -13.33 &  -26.06 &  -40.57 &  0.00827 &  0.02932 &  1.84 \\
$\mathsf{C2}$    &  1.17 &    1.90 &    0.91 &   -4.11 &   -9.87 &  -39.45 &  0.00155 &  0.32562 &  1.59 \\
$\mathsf{C3}$    &  1.21 &    1.02 &   -1.06 &   -8.18 &   -9.87 &  -23.91 &  0.01008 &  0.08807 &  1.53 \\
$\mathsf{C4}$    &  1.21 &    1.84 &   -1.02 &   -6.77 &   -9.87 &  -22.38 &  0.00230 &  0.03980 &  1.61 \\
$\mathsf{D1}$    &  1.08 &   41.75 &   -9.86 &  -14.06 &  -16.43 &  -39.39 &  0.00110 &  0.00526 &  1.71 \\
$\mathsf{D2}$    &  1.00 &   -8.52 &   -9.86 &  -16.43 &  -20.71 &  -39.40 &  0.00012 &  0.00967 &  1.71 \\
$\mathsf{D3}$    &  1.00 &   -1.99 &   -9.86 &  -16.44 &  -16.68 &  -39.39 &  0.00093 &  0.02219 &  1.71 \\
$\mathsf{D4}$    &  1.02 &   17.12 &   -9.86 &  -13.38 &  -16.40 &  -39.39 &  0.00015 &  0.02364 &  1.68 \\
$\mathsf{D5}$    &  1.06 &   35.55 &   -9.86 &  -13.60 &  -16.43 &  -38.18 &  0.00003 &  0.00642 &  1.69 \\
$\mathsf{D6}$    &  1.09 &   48.39 &   -9.86 &  -13.02 &  -15.31 &  -39.40 &  0.00097 &  0.01294 &  1.71 \\
$\mathsf{D7}$    &  1.27 &  305.11 &   -9.86 &  -11.95 &  -15.69 &  -37.38 &  0.00089 &  0.02138 &  1.74 \\
$\mathsf{D8}$    &  1.16 &   47.51 &   47.21 &   -9.86 &  -11.93 &  -27.74 &  0.00002 &  0.00433 &  1.68 \\
$\mathsf{D9}$    &  1.16 &   84.85 &    7.23 &   -9.86 &  -26.09 &  -26.26 &  0.00418 &  0.01027 &  1.68 \\
$\mathsf{D10}$   &  1.45 &  160.60 &   -9.86 &  -13.42 &  -16.15 &  -17.35 &  0.00001 &  0.01819 &  1.70 \\
$\mathsf{D11}$   &  1.37 &   61.75 &   -9.86 &  -15.43 &  -18.83 &  -37.26 &  0.00014 &  0.00139 &  1.77 \\
$\mathsf{D12}$   &  1.39 &   73.01 &   22.01 &   -9.86 &  -12.47 &  -20.32 &  0.00009 &  0.01913 &  1.74 \\
$\mathsf{E9a}$   &  1.06 &  239.23 &   67.69 &   32.76 &   -9.87 &  -10.11 &  0.00707 &  0.05096 &  1.60 \\
$\mathsf{E9b}$   &  2.92 &  495.15 &  290.82 &  238.04 &  178.01 &   86.36 &  0.00016 &  0.01846 &  1.53 \\
    \bottomrule
    \end{tabular}
    \caption{
    Results of the eigenvalue method for different scenarios from \cite{vg14}. Maximum errors $ \epsilon = \max_i \frac{\lVert S x_i - \lambda_i x_i \lVert}{\lambda_i x_i} $ and $ \epsilon_g = \max_i \frac{\left| \lambda_i - \lambda_{g,i} \right|}{\lambda_i} $ are overall smaller than $0.1$ (with the exception of scenario $\mathsf{C2}$). 
    }
    \label{table:eig-vg14}
\end{table}

\begin{figure}
\centering
\includegraphics[width=\textwidth]{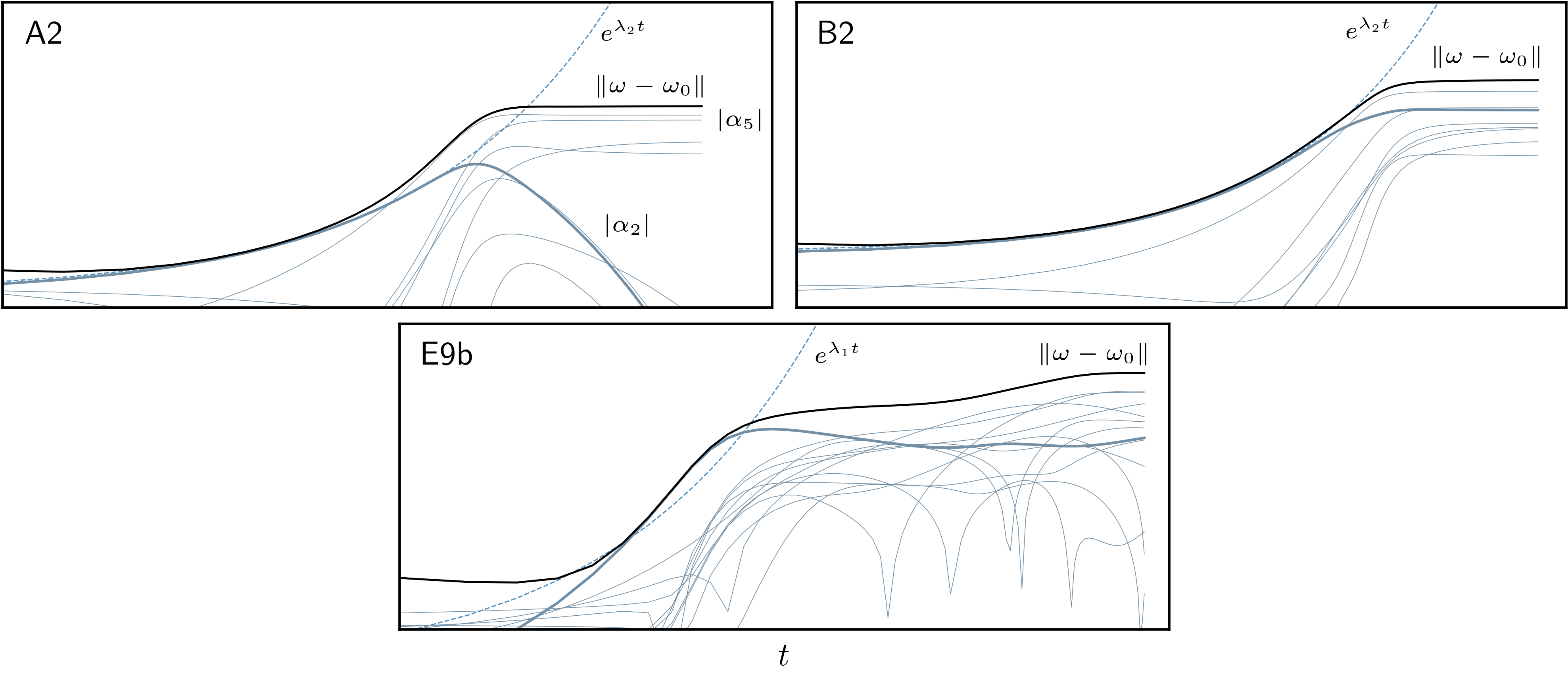}
\caption{
The norm of the perturbation $ \Vert \omega(t) - \omega_0 \Vert $ and the magnitude of the projections 
$ \lvert \alpha_j(t) \rvert $ are plotted against time. 
In scenario $\mathsf{E9b}$ the eigenvalue estimate of the growth rate is off by about $30 \%$  from the actual growth rate.
}
\label{fig:proiezioni}
\end{figure}

As for the connection between the spectrum and the evolution dynamics at early times, scenarios $\mathsf{A2}$, $\mathsf{B2}$ and $\mathsf{E9b}$ have been analyzed, after checking the (approximate) orthogonality of the eigenfunctions.
In figure \ref{fig:proiezioni}, the comparison between the actual time evolutions and the one predicted by the spectrum as $\sum_{i=0}^N \exp(\lambda_i t) e_i$ are represented. Here, we have used $N=8$ for cases $\mathsf{A2}$, $\mathsf{B2}$ and $N=12$ for $\mathsf{E9b}$.

\begin{figure}
\centering
\includegraphics{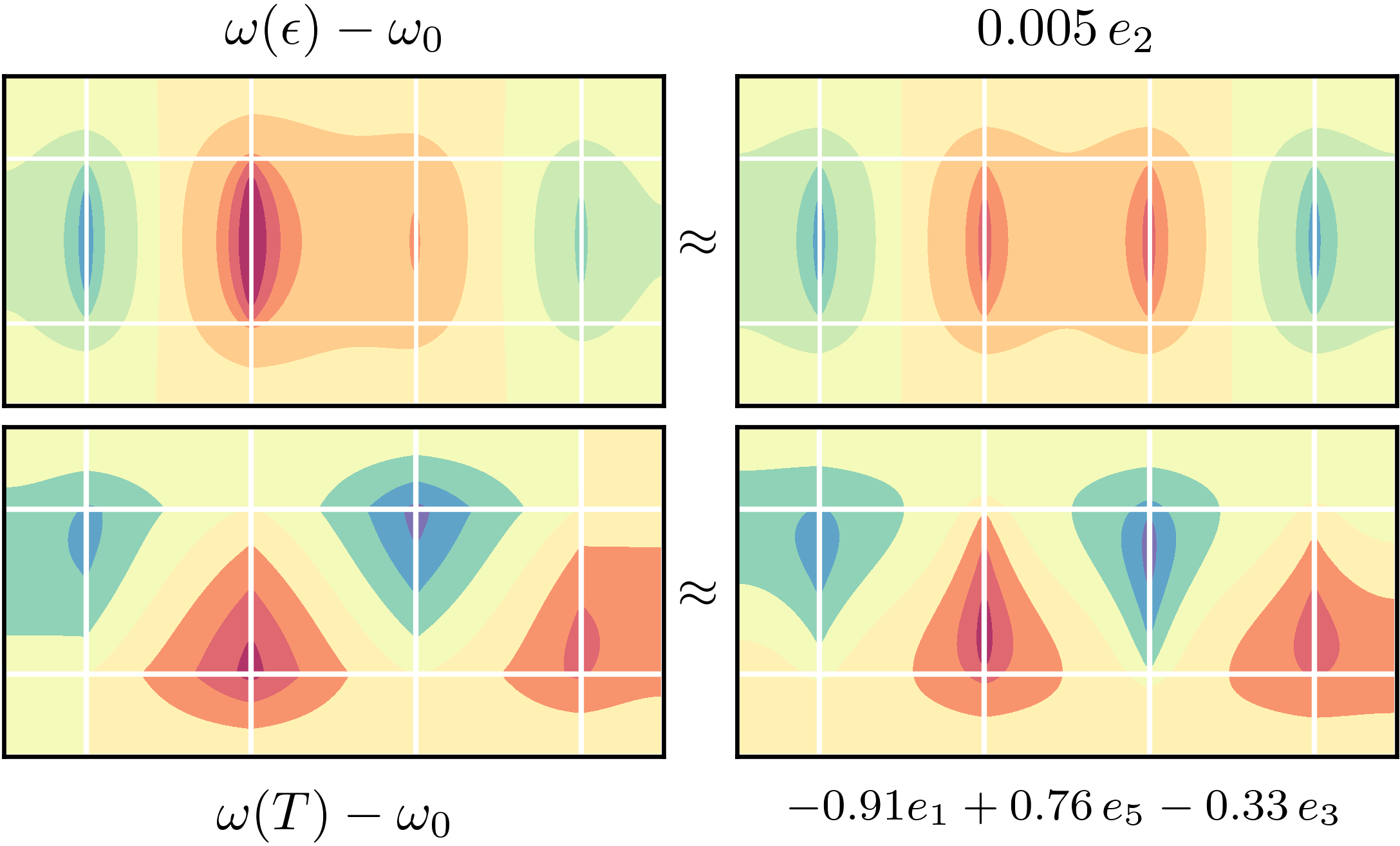}
\caption{At time $\epsilon$ (near zero), the growing perturbation has the form of one particular eigenfunction: 
$\delta\omega_\epsilon \approx \operatorname{proj}_{\langle e_2 \rangle} \delta\omega_\epsilon = \alpha_2 e_2$. The eigenvalue method however cannot predict the steady state profile, which is mostly associated to modes 2,3 and 5: $ \delta\omega_\infty \approx \operatorname{proj}_{\langle e_1,e_3,e_5 \rangle} \delta \omega_\infty $. 
}
\label{fig:proiezioni_n}
\end{figure}

The presence of multiple positive eigenvalues precludes the possibility of pinning down uniquely the shape and rate of growth of the instability. In all cases however, the growing perturbation has the form of one particular eigenfunction and its rate of growth is given by the associated eigenvalue, see figure \ref{fig:proiezioni_n}.
The eigenvalue with largest real part thus also provides a numerical upper bound on the rate of growth.

% \FloatBarrier
\subsection{The role of fracture circuits}
\label{sec:condensatore}

In \cite{vg14} the onset and strength of convection is linked to the existence of continuous fracture circuits, and different numerical experiments confirm its importance.
In this section we argue that quasi-continuous fracture circuits can also enable convective motions, for a surrounding porous matrix of sufficiently high permeability.

\begin{figure}
\centering
\includegraphics[width=0.8\textwidth]{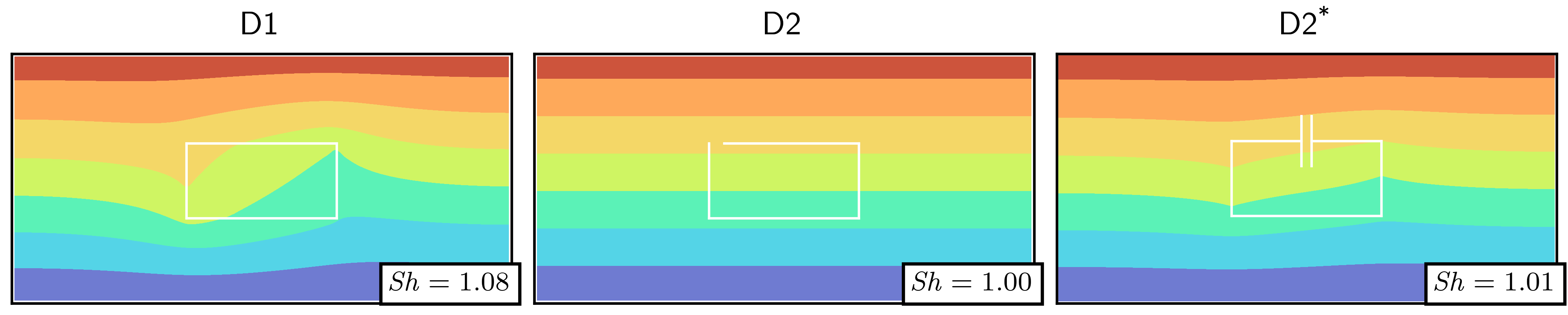}
\caption{ Cases $\mathsf{D1}$ and $\mathsf{D2}$ are compared in \cite{vg14} to show the necessity of fracture circuit continuity. Modifying the gap geometry as in $\mathsf{D2^*}$ restores the convective motion. } 
\label{fig:d2stella}
\end{figure}

Cases $\mathsf{D1}$ and $\mathsf{D2}$ in \cite{vg14} are compared to show the necessity of circuit continuity for convection to occur.
However, slightly changing the fracture geometry (parameters remain unchanged), as in scenario $\mathsf{D2^*}$ (see figure \ref{fig:d2stella}), we see that, although the strength of convection is significantly reduced, it is still possible.
The permeability of the surrounding medium is indeed large enough for the medium to be part of the convective circuit, 
though not large enough for convection to occur without the presence of fractures.

\begin{figure}
\centering
\includegraphics{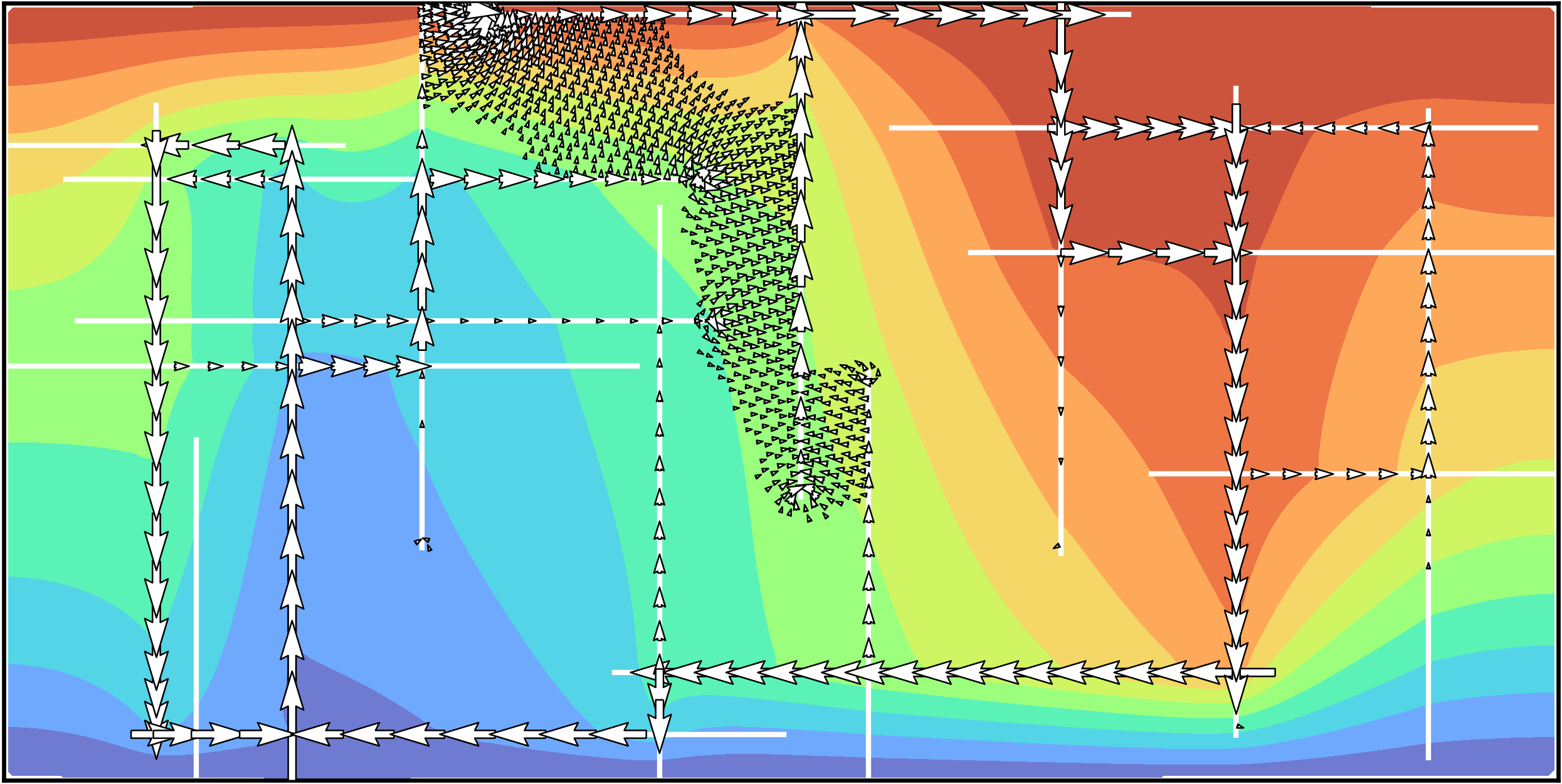}
\caption{Velocity $\bm{u}$ and concentration $\omega$ at steady state (scenario $\mathsf{E9b}$). The velocity vectors clearly indicate that the convective cell is not limited to the fracture network.}
\label{fig:e9b-vel}
\end{figure}

Scenario $\mathsf{D2^*}$ is purposefully built to highlight this phenomenon, we can however find similar configurations in more complex scenarios already introduced in \cite{vg14}.
Let us consider for instance network $\mathsf{E9b}$, shown in figure \ref{fig:e9b-vel}: velocity vectors plotted over the concentration profile clearly indicate, unlike what is claimed in \cite{vg14}, one large convective cell instead of two, separate ones. 
As in scenario $\mathsf{D2^*}$, the convective loop is not limited to the fracture network but crosses the porous matrix as well.

\begin{figure}
    \centering
    \includegraphics[width=\textwidth]{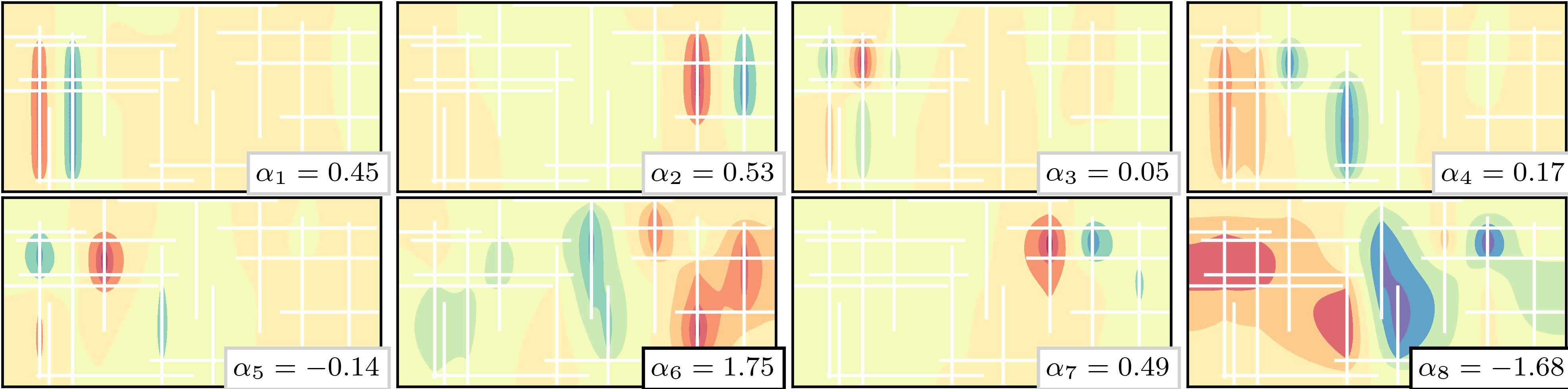}
    \caption{The first 8 modes and the corresponding projections $\alpha_i$ as defined by (\ref{eq:alpha}). Some modes describe convection around fracture circuites, such as $e_1$, $e_2$, $e_4$, $e_7$. At steady state however, dominant modes $e_6$ and $e_8$ both involve convective motion across the gaps.  } 
    \label{fig:e9b-proiezione}
\end{figure}

As further confirmation, we can decompose the solution of $\mathsf{E9b}$ using the approximate eigenvector basis (following the apprach outlined in section \ref{sec:exp-hrl}).
In figure \ref{fig:e9b-proiezione} we see how the dominant modes $e_6, e_8$, i.e. the modes corresponding to the largest magnitudes of $\alpha_i$, both involve fluid motion across the gap in the fracture network, through the porous matrix,
while modes looping around large continuous fracture circuits $e_1, e_2, e_4, e_7$ contribute only marginally to the steady state solution.

\begin{figure}    
    \centering
    \includegraphics[width=0.7\textwidth]{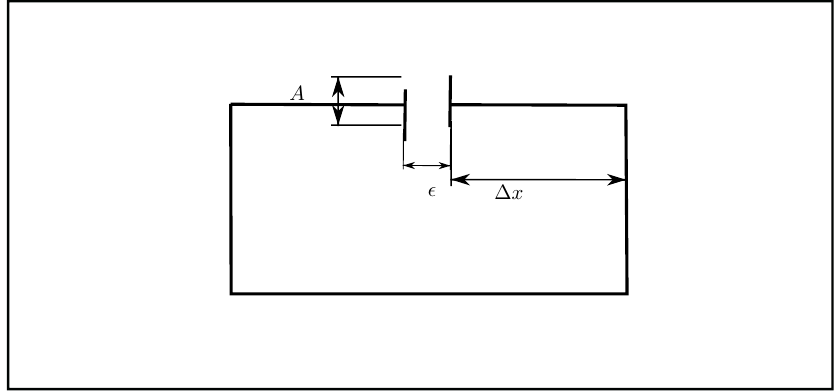}
    \caption{Discontinuous loop geometry}
    \label{fig:condensatore-geo}
\end{figure}

\begin{figure}
    \centering
    \includegraphics[width=\textwidth]{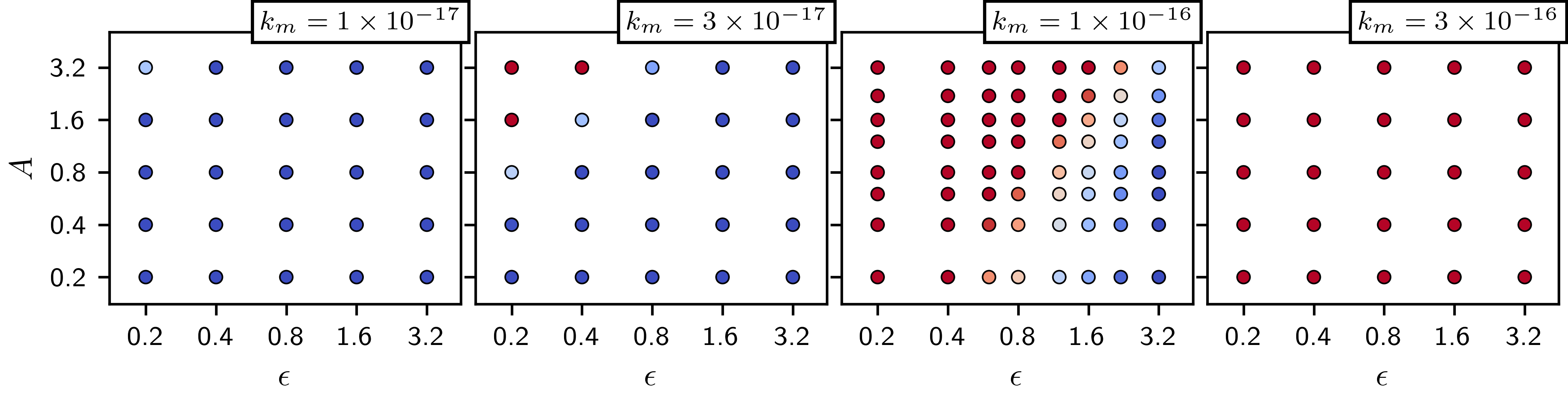}
    \caption{Possibility of convection (positive largest eigenvalue) for varying $\epsilon, A, k_m$ is indicated by a red dot. Meaning of geometric parameters $\epsilon, A$ is indicated in figure \ref{fig:condensatore-geo}. }
    \label{fig:condensatori-val}
\end{figure}

Let us parametrize the geometry of a fracture circuit with a gap as in figure \ref{fig:condensatore-geo}, to relate it to the possibility of convection with the analysis of eigenvalues sign. Results are reported in figure \ref{fig:condensatori-val}.
As expected, for very small matrix permeability gaps stop flow, thus inhibiting convection.
On the other hand, for larger matrix permeability $ k_m = \num{3e-16} $, the corresponding Rayleigh number is $ \Ra \approx 20 $, very near the critical $ \Ra_c = 4 \pi^2 $ for which the matrix can exhibit convection even without the aid of fractures; convection is thus possible across very large gaps. 
For matrix permeability between these two extremes, the geometry of the gap determines whether convection across is possible or not.

This behavior can be explained, at least qualitatively, focusing on the top (broken) edge of the fracture circuit composed by the two fracture segments of length $\Delta x$ (segments $A$ and $C$), and the gap of width $\epsilon$ (denoted by $B$). Assuming that across the gap flow occurs only across an area $A$, and that exchange with the porous matrix elsewhere re negligible, to ensure convection in the circuit we have to guarantee that the conductivity $G_B\geq G_A=G_C$, which translates into

$$\dfrac{A k_m}{\epsilon}\geq \dfrac{ k_f b}{\Delta x},$$

or, equivalently, 

$$\dfrac{k_m A}{k_f b}\dfrac{\Delta x}{\epsilon}\geq 1.$$

\begin{figure}
\centering
\includegraphics{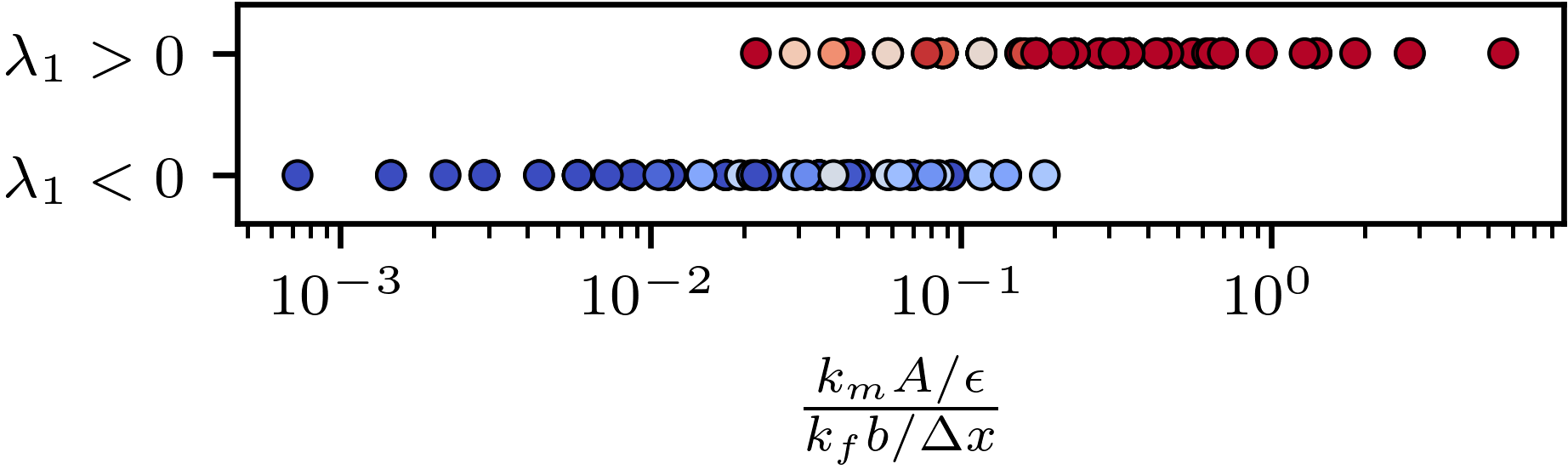}
\caption{
Each dot corresponds to a particular combination of $k_m, \epsilon, A$. 
As in figure \ref{fig:condensatori-val}, red and blue dots correspond respectively to a positive and negative largest eigenvalue. 
}
\label{fig:condensatori-res}
\end{figure}

This relation correctly predicts the qualitative behaviour of the numerical simulations reported in figure \ref{fig:condensatori-val}: small gaps with large surface areas favour convection across the gap. 
The results of a more quantitative evaluation, checking the sign of eigenvalues for different combinations of the parameters, are illustrated in figure \ref{fig:condensatori-res}. Although the threshold does not appear to be so sharp, and its numerical value is closer to $10^-1$ than to $1$, the horizontal separation of positive leading eigenvalues from negative ones indicates that the model is capturing part of the physics of the phenomenon.

\subsection{Three-dimensional HRL problem}
\label{sec:exp-hrl-3d}

Both the direct method (exposed in section \ref{sec:direct}) and the eigenvalue method (section \ref{sec:eig}) make no assumptions on the dimensionality of the ambient space.
In this section, we are going to apply the latter to simple three-dimensional scenarios.

\begin{figure}
    \centering
\includegraphics[width=\textwidth]{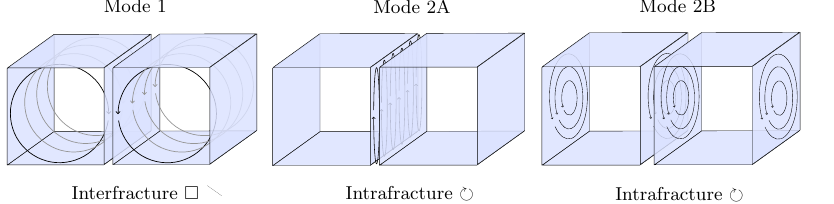}
\caption{Classification of convective modes presented in \cite{nield}. Two-dimensional simulations can only capture the interfracture mode 1. Three-dimensional simulations can also capture intrafracture modes 2A and 2B, unless fractures are treated as lower-dimensional subdomains: in this case mode 2A is impossible to capture.}
    \label{fig:3d-modes}
\end{figure}

In \cite{nield}, a classification of possible modes of convection in fractured porous media is proposed.
The paper further explains how the assumption of a two-dimensional domain severely limits the applicability to experimental cases.
This is due to the three-dimensional character of the dominant mode of convection -- mode 2B in figure \ref{fig:3d-modes} -- which cannot be captured by a two-dimensional analysis.

In \cite{vg15}, different three-dimensional scenarios are analyzed and numerically simulated. 
The numerical results from the regular three-dimensional fracture circuit (here denoted as scenario 6), show that when convection is available (for large enough fracture aperture), the dominant mode is the interfracture mode 1.

The geometry simulated in \cite{vg15} being practically identical to the geometry studied in \cite{nield}, the two results stand in contradiction.
Indeed, \cite{vg15} acknowledges the contradiction, 
attributing it to both the matrix-fracture coupling conditions and to the Rayleigh averaging strategy used in the analysis of \cite{nield}.
In particular, the use of a Rayleigh number based on an averaged permeability was shown to be ineffective in the case of low density fracture networks in \cite{vg14} (as briefly discussed in section \ref{sec:frac-hrl}).

We analyze the 3D problem with the eigenvalue method for different scenarios presented in \cite{vg15}, and collect the results in table \ref{table:vg15}, where we report for different geometries and different fracture apertures the unstable modes, and compare them with the convection modes predicted by \cite{vg15}. 

As in the two-dimensional case, the instability thresholds are identical to the ones obtained in \cite{vg15}, apart from a slight misprediction for scenario $\mathsf{9a}$.

In all scenarios furthermore, the eigenfunctions clearly follow the classification proposed in \cite{nield}, corresponding to either interfracture $\square$ or intrafracture $\circlearrowright$ modes.

\begin{table}
\small
\centering 
\begin{tabular}{ cccc }\toprule 
Case & $b \times \SI[print-unity-mantissa=false]{1e-5}{\meter}$ & Mode \cite{vg15} & Eigenfunctions  \\
\midrule
$\mathsf{5}$   & $ 1.8\,(1-\epsilon) $ & ---                              & ---                                                              \\
$\mathsf{5}$   & $ 1.8\,(1+\epsilon) $ & $\circlearrowright$              & $2\times\bm{\circlearrowright}$                                  \\
$\mathsf{5}$   & $ 2.0               $ & $\circlearrowright$              & $4\times\bm{\circlearrowright}$                                  \\
$\mathsf{6}$   & $ 1.4\,(1-\epsilon) $ & ---                              & ---                                                              \\
$\mathsf{6}$   & $ 1.4\,(1+\epsilon) $ & $\square$                        & $1\times\bm{\square} \quad 1\times\circlearrowright$             \\
$\mathsf{6}$   & $ 1.6               $ & $\square$                        & $1\times\square \quad 2\times\bm{\circlearrowright}$             \\
$\mathsf{6}$   & $ 1.8               $ & $\square$                        & $1\times\square \quad 4\times\bm{\circlearrowright}$             \\
$\mathsf{6}$   & $ 2.0               $ & $\square$                        & $1\times\square \hspace{5pt}  12\times\bm{\circlearrowright}$    \\
$\mathsf{9a}$  & $ 1.5               $ & ---                              & ---                                                              \\
$\mathsf{9a}$  & $ 1.6               $ & ---                              & ---                                                              \\
$\mathsf{9a}$  & $ 1.7               $ & ---                              & $ 1 \times \bm{\circlearrowright} $                              \\
$\mathsf{9a}$  & $ 1.8               $ & $ \square \; \circlearrowright $ & $ 1\times\square \quad 2\times\bm{\circlearrowright} $           \\
$\mathsf{9a}$  & $ 1.9               $ & $ \square \; \circlearrowright $ & $ 1\times\square \quad 2\times\bm{\circlearrowright} $           \\
$\mathsf{9a}$  & $ 2.0               $ & $ \square \; \circlearrowright $ & $ 1\times\square \quad 5\times\bm{\circlearrowright} $           \\
$\mathsf{9b}$  & $ 1.5               $ & ---                              & ---                                                              \\
$\mathsf{9b}$  & $ 1.6               $ & ---                              & ---                                                              \\
$\mathsf{9b}$  & $ 1.7               $ & $ \square $                      & $ 1\times\square \quad 1\times\bm{\circlearrowright} $           \\
$\mathsf{9b}$  & $ 1.8               $ & $ \square $                      & $ 1\times\square \quad 2\times\bm{\circlearrowright} $           \\
$\mathsf{9b}$  & $ 1.9               $ & $ \square $                      & $ 1\times\square \quad 2\times\bm{\circlearrowright} $           \\
$\mathsf{9b}$  & $ 2.0               $ & $ \square $                      & $ 1\times\square \quad 6\times\bm{\circlearrowright} $           \\
\bottomrule\end{tabular}    
\caption{$ \epsilon = 0.05 $. In scenarios $\mathsf{5}$, $\mathsf{6}$ and $\mathsf{9b}$ the predicted instability threshold agrees with \cite{vg15}. The threshold is slightly ($\lambda = 1.15$) mispredicted for scenario $\mathsf{9a}$. Between interfracture $\square$ and intrafracture $\circlearrowright$ modes, the dominant one is printed in bold. $ \, \square \, \circlearrowright $ indicates a combination of interfracture and intrafracture modes. }
\label{table:vg15}
\end{table}

In particular, for scenario 6 we observe that (\romannumeral 1) the interfracture mode is first unstable one appearing at $ b \approx \SI{1.4e-5}{\meter} $
(\romannumeral 2) already at $ b \approx \SI{1.6e-5}{\meter} $, different intrafracture modes are available and dominate over the interfracture mode:
$ \lambda^\circlearrowright_1 = 21.0,\,\lambda^\circlearrowright_2 = 20.1,\,\lambda^\square = 17.6 $
(\romannumeral 3) intrafracture convection modes remain dominant for larger fracture apertures.

Thus, excluding apertures for which the setup is very near to the instability threshold, intrafracture convection modes dominate over interfracture convection modes when both modes are available.
This conclusion also holds for scenarios $\mathsf{9a}$ and $\mathsf{9b}$.

Although our results seem to confirm the results of \cite{nield}, where the intrafracture mode is also identified as dominant,
they do not contradict the results in \cite{vg15}. 
Indeed the eigenvalue analysis is strictly localized around the diffusive equilibrium solution. In particular, the eigenvalues give no indication on which modes will be present at steady state. 
As seen in e.g. scenario $\mathsf{A2}$ (figure \ref{fig:proiezioni_n}), even modes with decaying behaviour near equilibrium may turn out to be dominant at steady state.

\begin{figure}
    \centering
    \includegraphics[width=\textwidth]{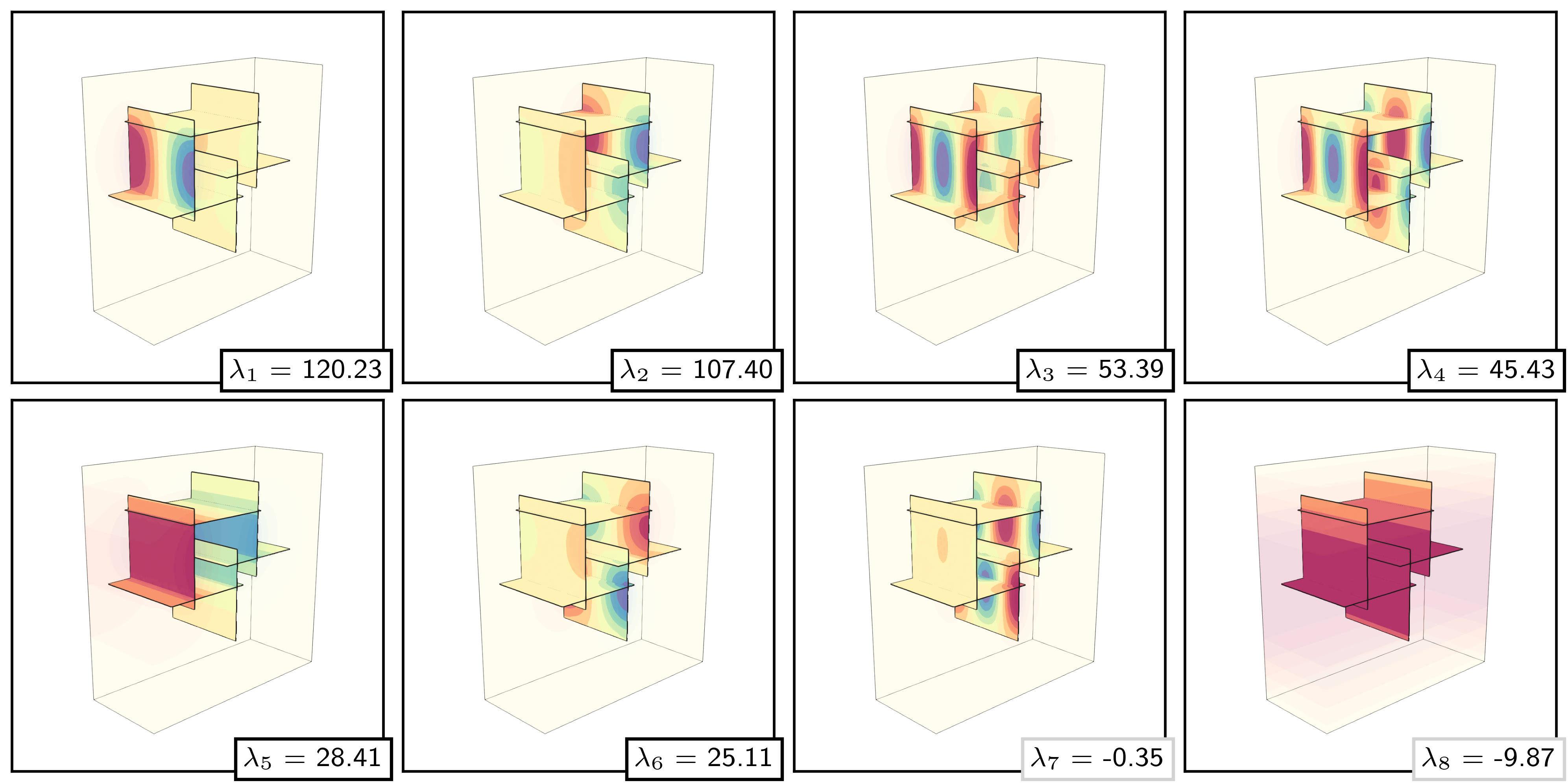}
    \caption{Results of the eigenanalysis applied to scenario $\mathsf{9a} \; (b = 2.0)$. The eigenfunction corresponding to the largest eigenvalue (referred to as dominant mode) is an intrafracture mode. Eigenfunction $e_5$ is the only interfracture mode. }
    \label{fig:vg15_9aF}
\end{figure}

% \FloatBarrier
\subsection{Computational considerations}
\label{sec:computational}

We want to discuss some of the computational aspects of the eigenvalue method, comparing it in particular to the direct method of assessing stability in terms of computational cost.

The computationally expensive steps of the eigenvalue method can be divided in 
(\emph{a}) grid construction, initial operator discretization and Newton iterations to reach steady state (needed as initial condition)
(\emph{b}) construction of the S operator defined in (\ref{eq:eig}) (including LU factorization of $A_{yy}$)
(\emph{c}) computation of $k$ largest eigenvalues.
Step (\emph{c}) can itself be subdivided in the two costly operations (\emph{i}) matrix vector products $S v$ and (\emph{ii}) orthogonalization.
We will initially compare the computational cost of assessing the possibility of convection using the eigenvalue method against the direct method.
Next, we will discuss how the cost changes with number of degrees of freedom and number of sought eigenvalues.

In section \ref{sec:exp-hrl}, we saw how the direct method and the eigenvalue method can be used in a complementary fashion,
each of them giving answers inaccessible to the other.
We can however restrict the question to the possibility of convection since both methods are a valid way of finding an answer (e.g. see table \ref{table:eig-vg14}).
For a fairer comparison, the direct forward simulations will be as soon as the Sherwood number exceeds 1, which would already indicate the possibility of convection.
Table \ref{table:prof-comp} reports the timing comparison and the relative time spent in each of the most costly operations for different test cases.

The cost of the direct method 
is dominated by the assembly of the problem Jacobian in the two-dimensional cases
and by the linear solves in the three-dimensional cases (which possess a much larger number of degrees of freedom).
Both of these operations are performed for each Newton iteration. 
The number of Newton iterations per timestep advancement mostly varies between 2 and 3 for all the simulations.
Relaxing the convergence tolerance of the Newton method is one of the options to speed up the method, at the cost of lower accuracy of the solution.
This might be a valid option considering that what we are interested in is whether convection is possible, not in the detailed behaviour of the solution.

The eigenvalue method is much more efficient with respect to the direct method in low fracture density cases.
For test cases with complex fracture configurations such as $\mathsf{E9a}$ and $\mathsf{E9b}$, the performance of the two methods (direct and eigenvalue) approximately matches.

In the different test cases, the eigenvalue method dedicates approximately equal times to orthogonalization and matrix-vector products.
Their relative weights can actually be adjusted by varying the size $m$ of the Krylov basis in the Krylov-Schur algorithm:
for every iteration of the algorithm, the cost of orthogonalization scales with $m^2$ while the cost of the matrix-vector products scales with $m$.
Tuning $m$ for the problem at hand may be a possible option for improving the efficiency of the method.

\begin{table}\footnotesize\centering
\begin{tabular}{lr|rccc|rcccc}\toprule
Scenario       & $\mathit{NoE}$  & direct (s) & \texttt{init} & \texttt{assembly}  & \texttt{linsolve}   & eig (s) & \texttt{init} & \texttt{lu} & \texttt{matvec} & \texttt{ortho} \\ \midrule
$\mathsf{A1}$  &  8192           &     127.1  & 0.03  & 0.70    & 0.15   &   22.5 & 0.17 & 0.02 & 0.13 & 0.23  \\
$\mathsf{A2}$  &  8192           &     186.7  & 0.02  & 0.74    & 0.16   &   23.0 & 0.24 & 0.03 & 0.11 & 0.18  \\
$\mathsf{B1}$  & 10016           &     249.4  & 0.03  & 0.73    & 0.17   &  37.0 & 0.17 & 0.02 & 0.29 & 0.28  \\
$\mathsf{B2}$  & 10422           &     260.0  & 0.03  & 0.72    & 0.20   &   53.0 & 0.17 & 0.02 & 0.33 & 0.34  \\
$\mathsf{D1}$  &  5000           &     125.0  & 0.01  & 0.77    & 0.13   &   20.0 & 0.22 & 0.03 & 0.05 & 0.08  \\
$\mathsf{D2}$  &  5000           &      87.9  & 0.03  & 0.75    & 0.12   &   13.9 & 0.21 & 0.02 & 0.09 & 0.15  \\
$\mathsf{E9a}$ & 11379           &     160.1  & 0.08  & 0.68    & 0.18   &  165.5 & 0.07 & 0.00 & 0.40 & 0.42  \\
$\mathsf{E9b}$ & 12073           &     137.3  & 0.11  & 0.59    & 0.15   &  107.5 & 0.21 & 0.02 & 0.36 & 0.28  \\
$\mathsf{6A}$  & 38400           &    1053.9  & 0.02  & 0.21    & 0.74   &  172.2 & 0.18 & 0.11 & 0.52 & 0.13  \\
$\mathsf{6C}$  & 38400           &     886.1  & 0.04  & 0.22    & 0.72   &  191.7 & 0.22 & 0.12 & 0.48 & 0.12  \\
\bottomrule
\end{tabular}
\caption{Comparison of time to assess the possibility of convection through direct method and eigenvalue method. The relative times of their more costly operations are also reported. }
\label{table:prof-comp}
\end{table}

\begin{figure}
\centering
\includegraphics{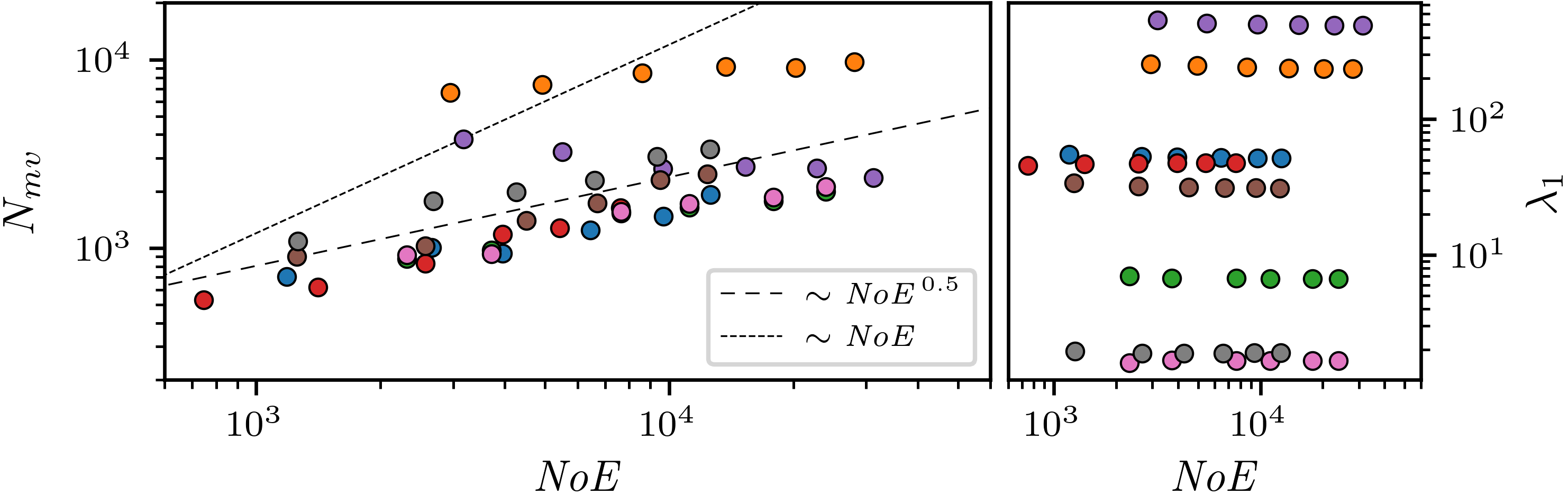}
\caption{On the left, number of matrix-vector products $N_{mv}$ plotted against number of elements $\noe$ on the left. 
Each different color corresponds to a different test case among $ \{ \mathsf{A2}, \mathsf{B2}, \mathsf{C2}, \mathsf{D8}, \mathsf{E9a}, \mathsf{E9b}, \mathsf{6A}, \mathsf{6C} \} $. On the right, we check that eigenvalues are converging (or have converged) for increasing grid resolution. }
\label{fig:prof-eig-2}
\end{figure}

In figure \ref{fig:prof-eig-2} we illustrate how the number of matrix-vector products $N_{mv}$ changes with the number of mesh elements $\noe$.
The results show that $N_{mv}$ depends weakly ($N_{mv} \sim \noe^{0.5} $) on the size of the problem (although a lot of variability remains).
Note that however the cost of each matrix-vector product $S v$ scales as $C_{mv} \sim \noe^2$: the cost of each matrix-vector product $S v$ scales as $C_{mv} \sim N_{mv}^2$: the presence of $A_{yy}^{-1}$ in (\ref{eq:eig}) makes the $S$ matrix dense.
The scaling estimate for the total cost of the matrix-vector products $\mathit{TC}_{mv}$ is thus
$$ \mathit{TC}_{mv} = N_{mv} C_{mv} \sim \noe^{0.5} \noe^2 = \noe^{2.5} . $$
The estimate is consistent with the data reported in table \ref{table:prof-comp}, for which a least-squares estimate gives $ \mathit{TC}_{mv} \sim \noe^{2.2 \pm 0.4} $.

\begin{figure}
    \centering
    \includegraphics{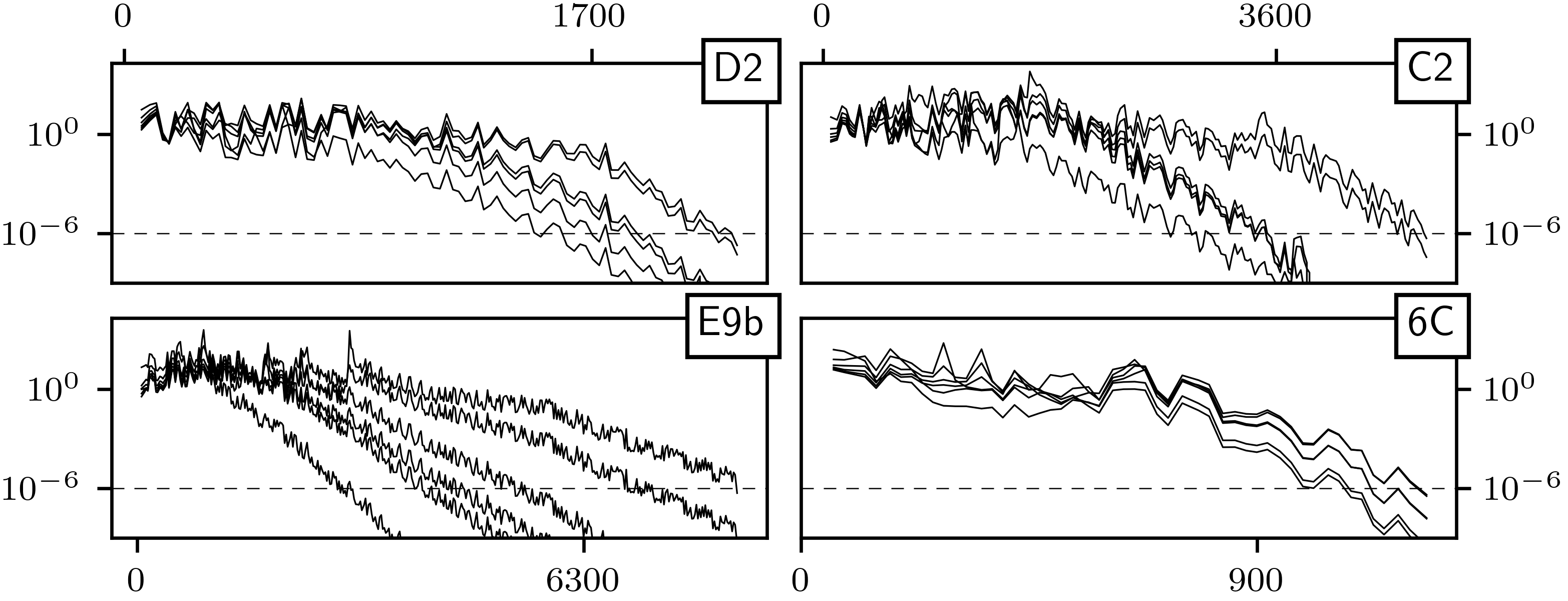}
    \caption{ 
Each line corresponds to the relative error for an eigenvalue estimate plotted against number of matrix-vector products $N_{mv}$.
}
    \label{fig:prof-eig-1}
\end{figure}

Finally, the convergence of the  Krylov-Schur algorithms is represented for four different test cases during the earch for multiple eigenvalues. As shown in figure \ref{fig:prof-eig-1} in all the test cases, once the algorithms pins down an eigenvalue, the corresponding error starts decreasing at exponential rate. Depending on the test case, we can have scenarios in which the eigenvalues all practically converge together, as opposed to cases where computing each successive eigenvalue requires considerably more computation.

\section{Conclusion}
\label{sec:conclusioni}

The aim of this study was to better understand the influence of fractures on the possibility of free convection in porous media.
To this aim, we have described a mathematical model for density driven flow in the presence of fractures,
and the corresponding numerical approximation. In addition to the direct ''forward'' numerical solution of the problem
we have proposed and implemented a novel method for the assessment of convective stability through the eigenvalue analysis of the linearized numerical problem.

The new method is shown to be in agreement with existing literature cases both in simple and complex fracture configurations.
With respect to direct simulation in time, its results provide further information on the possibility of convection.
In particular, as shown in section \ref{sec:exp-hrl} the sign of the leading eigenvalue correctly predicts the onset of convection, and its magnitude provides quantitative estimates of the rate of growth or decay of perturbations.
Furthermore the computational cost of the method has proven to be in the worst cases equal, and in the best cases up to an order of magnitude faster than the direct solution method.

The fact that the eigenvalue method closely mimics the analytical method of investigating stability clarifies what inferences can be made from the the results of a linear stability analysis. In particular, analyzing stability through the study of the system around equilibrium solutions (both numerically and analytically) is restrictive in that it cannot predict whether, far from the equilibrium solution, transient convective motion is possible and which convective modes will be dominant at steady state.

Moreover, an in depth study of the particular scenario E9b in \cite{vg14} has further complicated the question of analyzing free convection in the presence of fractures: for realistic sets of problem parameters, the porous matrix is indeed able to participate in the convective motion. 
Thus, stability criteria based on the fracture network alone, e.g. the presence of large continuous fracture circuits as a trigger for convection, are shown to be somewhat limited in applicability.

Given the results of this study, we could expand the work in different directions.

On the numeric side, the method used to compute the eigenvalues at the moment does not take into account the structure of the particular problem.
However, the $S$ matrix is built up from submatrices of the Jacobian associated to the linearized problem.
Therefore, a close study of the Jacobian may suggest ways of speeding up the computation of the eigenvalues by exploiting the structure of the $S$ matrix.

Finally, though the eigenvalue method has proven to be faster than the direct method of assessing convective stability, the method could be further sped up if we are willing to sacrifice accuracy, since, for this purpose, we are only interested in the sign of the most positive eigenvalue. If during the computation of the eigenvalues we find one to be positive, with error reasonably smaller than its magnitude, we can already predict the possibility of convective motion.
In cases where the method is applicable, such as the fractured HRL problem, the additional speed 
may enable a more systematic (or even statistical) study of how global properties of fracture networks such as fracture density, connectivity and characteristic sizes are related to the possibility of convective motion.

\section*{Acknowledgements}
The authors acknowledge the support by MUR, grant Dipartimento di Eccellenza 2023–2027.


\begin{thebibliography}{10}

\bibitem{aavatsmark02}
Ivar Aavatsmark.
\newblock An {Introduction} to {Multipoint} {Flux} {Approximations} for
  {Quadrilateral} {Grids}.
\newblock {\em Computational Geosciences}, 6(3/4):405--432, 2002.

\bibitem{diersch02}
H.-J.G. Diersch and O.~Kolditz.
\newblock Variable-density flow and transport in porous media: approaches and
  challenges.
\newblock {\em Advances in Water Resources}, 25(8-12):899--944, August 2002.

\bibitem{elder67}
J.~W. Elder.
\newblock Steady free convection in a porous medium heated from below.
\newblock {\em Journal of Fluid Mechanics}, 27(1):29--48, January 1967.

\bibitem{hortonrogers}
C.~W. Horton and F.~T. Rogers.
\newblock Convection {Currents} in a {Porous} {Medium}.
\newblock {\em Journal of Applied Physics}, 16(6):367--370, June 1945.

\bibitem{porepy}
Eirik Keilegavlen, Runar Berge, Alessio Fumagalli, Michele Starnoni, Ivar
  Stefansson, Jhabriel Varela, and Inga Berre.
\newblock {PorePy}: an open-source software for simulation of multiphysics
  processes in fractured porous media.
\newblock {\em Computational Geosciences}, 25(1):243--265, February 2021.

\bibitem{lapwood_1948}
E.~R. Lapwood.
\newblock Convection of a fluid in a porous medium.
\newblock {\em Mathematical Proceedings of the Cambridge Philosophical
  Society}, 44(4):508--521, October 1948.

\bibitem{martin05}
Vincent Martin, Jérôme Jaffré, and Jean~E. Roberts.
\newblock Modeling {Fractures} and {Barriers} as {Interfaces} for {Flow} in
  {Porous} {Media}.
\newblock {\em SIAM Journal on Scientific Computing}, 26(5):1667--1691, January
  2005.

\bibitem{rees}
D.A.S. Rees and L.~Storesletten.
\newblock The onset of convection in a two-layered porous medium with
  anisotropic permeability.
\newblock {\em Transport in Porous Media}, 128, 2019.

\bibitem{shafabakhsh19}
Paiman Shafabakhsh, Marwan Fahs, Behzad Ataie-Ashtiani, and Craig~T. Simmons.
\newblock Unstable {Density}-{Driven} {Flow} in {Fractured} {Porous} {Media}:
  {The} {Fractured} {Elder} {Problem}.
\newblock {\em Fluids}, 4(3):168, September 2019.

\bibitem{shikaze98}
Steven~G Shikaze, E.A Sudicky, and F.W Schwartz.
\newblock Density-dependent solute transport in discretely-fractured geologic
  media: is prediction possible?
\newblock {\em Journal of Contaminant Hydrology}, 34(3):273--291, October 1998.

\bibitem{nield}
C.~T. Simmons, J.~M. Sharp, and D.~A. Nield.
\newblock Modes of free convection in fractured low‐permeability media.
\newblock {\em Water Resources Research}, 44(3):2007WR006551, March 2008.

\bibitem{starnoni19}
M.~Starnoni, I.~Berre, E.~Keilegavlen, and J.~M. Nordbotten.
\newblock Consistent {MPFA} {Discretization} for {Flow} in the {Presence} of
  {Gravity}.
\newblock {\em Water Resources Research}, 55(12):10105--10118, December 2019.

\bibitem{stathopoulos98}
Andreas Stathopoulos, Yousef Saad, and Kesheng Wu.
\newblock Dynamic {Thick} {Restarting} of the {Davidson}, and the {Implicitly}
  {Restarted} {Arnoldi} {Methods}.
\newblock {\em SIAM Journal on Scientific Computing}, 19(1):227--245, January
  1998.

\bibitem{stewart02}
G.~W. Stewart.
\newblock A {Krylov}--{Schur} {Algorithm} for {Large} {Eigenproblems}.
\newblock {\em SIAM Journal on Matrix Analysis and Applications},
  23(3):601--614, January 2002.

\bibitem{voss1987variable}
Clifford~I. Voss and William~R. Souza.
\newblock Variable density flow and solute transport simulation of regional
  aquifers containing a narrow freshwater‐saltwater transition zone.
\newblock {\em Water Resources Research}, 23(10):1851--1866, October 1987.

\bibitem{vg15}
Katharina Vujević and Thomas Graf.
\newblock Combined inter- and intra-fracture free convection in fracture
  networks embedded in a low-permeability matrix.
\newblock {\em Advances in Water Resources}, 84:52--63, October 2015.

\bibitem{vg14}
Katharina Vujević, Thomas Graf, Craig~T. Simmons, and Adrian~D. Werner.
\newblock Impact of fracture network geometry on free convective flow patterns.
\newblock {\em Advances in Water Resources}, 71:65--80, September 2014.

\end{thebibliography}
\end{document}